\documentclass[10pt]{article}

\title{\textbf{Superconvergent Non-Polynomial Approximations}
\thanks{Last updated on June 11, 2021.
The research of the authors was supported in part by AFOSR grants FA9550-19-1-0281 and
FA9550-17-1-0394 and NSF grant DMS 1912183.}
}

\author{Andrew J. Christlieb \thanks{Department of Computational Mathematics, Science and Engineering, Michigan State University, East Lansing, MI, 48824, United States; \href{mailto:christli@msu.edu}{christli@msu.edu}.}
\and William A. Sands \thanks{Department of Computational Mathematics, Science and Engineering, Michigan State University, East Lansing, MI, 48824, United States; \href{mailto:sandswi3@msu.edu}{sandswi3@msu.edu}.}
\and Hyoseon Yang  \thanks{Department of Computational Mathematics, Science and Engineering, Michigan State University, East Lansing, MI, 48824, United States; \href{mailto:hyoseon@msu.edu}{hyoseon@msu.edu} (corresponding author).}
}

\date{}

\usepackage{amsmath}
\usepackage{epsfig}
\usepackage{color,graphics}
\usepackage{booktabs}
\usepackage{epstopdf}
\usepackage[all]{xy}
\usepackage{fixmath, bm, amssymb, amsfonts, latexsym,multirow}
\usepackage[a4paper,margin=1in]{geometry}
\usepackage{upgreek,comment,setspace,subfigure}
\usepackage[toc,page]{appendix}
\usepackage{algorithm}
\usepackage[noend]{algpseudocode}
\usepackage{hyperref}
\hypersetup{
    colorlinks=true,
    linkcolor=blue,
    citecolor=blue,
    filecolor=blue,      
    urlcolor=blue,
}

\newtheorem{theorem}{Theorem}[section]

\newtheorem{definition}[theorem]{Definition}
\newtheorem{example}[theorem]{Example}
\newtheorem{remark}[theorem]{Remark}

\newcommand{\baa}{\begin{eqnarray*}}
	\newcommand{\eaa}{\end{eqnarray*}}
\newcommand{\ba}{\begin{equation}}
\newcommand{\ea}{\end{equation}}

\newcommand{\RR}{\mathbb R}
\newcommand{\NN} {\mathbb N}

\newcommand{\CC} {\mathbb C}

\newcommand{\erf}{\text{erf}}
\newcommand{\etal}{\textit{et al.}}

\begin{document}

    \maketitle
	
	\begin{abstract}
In this paper, we introduce a superconvergent approximation method that employs radial basis functions (RBFs) in the numerical solution of conservation laws.  The use of RBFs for interpolation and approximation is a well developed area of research. Of particular interest in this work is the development of high order finite volume (FV) weighted essentially non-oscillatory (WENO) methods, which utilize RBF approximations to obtain required data at cell interfaces. Superconvergence is addressed through an analysis of the truncation error, resulting in expressions for the shape parameters that lead to improvements in the accuracy of the approximations. This study seeks to address the practical elements of the approach, including the evaluations of shape parameters as well as  hybrid implementation. To highlight the effectiveness of the non-polynomial basis, in shock-capturing, the proposed methods are applied to one-dimensional hyperbolic and weakly hyperbolic systems of conservation laws and compared with several well-known FV WENO schemes in the literature. In the case of the non-smooth, weakly hyperbolic test problem, notable improvements are observed in predicting the location and height of the finite time blowup. The convergence results demonstrate that the proposed schemes attain notable improvements in accuracy, as indicated by the analysis of the reconstructions. We also include a discussion regarding extensions to higher dimensional problems, along with convergence results for a  nonlinear scalar problem.
\end{abstract}




	
\section{Introduction}
This work concerns the development of numerical schemes to solve conservation laws, which are of the form
\begin{align}\label{eq:claw0}
\begin{cases}
& u_t + \nabla \cdot F(u) = 0,\quad   \ x\in \RR^m,\\
&     u(x,0)    = u_0(x),
\end{cases}
\end{align}
where $F(u)$ is the flux function and $ u_0(x)$ is the prescribed initial data. The necessity and computational cost of high order numerical methods for the problem \eqref{eq:claw0} have motivated, to a large extent, the development of so-called \textit{superconvergent} approximation schemes. Superconvergent schemes can be broadly defined as those which converge at a rate faster than what is theoretically expected. Henceforth, we shall use this definition when referring to superconvergence for the approximations appearing in this work.  

Finite volume (FV) schemes are among the most popular methods used to solve hyperbolic conservation laws. An attractive feature of such schemes is that they evolve cell average data, which makes the discretization naturally conservative. These conservation properties also make them suitable for use in, for example, adaptive mesh refinement algorithms, where data needs to be frequently transferred between levels within a hierarchy of grids. The use of cell average data, rather than point-wise values, greatly simplifies the task of making such transfers conservative. While there exist many approaches to obtain high order FV discretizations, we restrict our focus to the class of so-called weighted ENO (WENO) schemes \cite{Jiang1, L2}, introduced by Jiang and Shu \cite{Jiang1}, which shall be referred to as classical WENO or WENO-JS. These methods were built on the success, and, in a sense, the limitations of the essentially non-oscillatory (ENO) schemes \cite{Har4, Har5, shu1, shu2}, which employed an adaptive interpolation stencil on a small set of available candidate points, to construct approximations at cell interfaces. WENO approaches, such as \cite{Jiang1}, in contrast, make use of all available candidate points in the reconstruction through a convex combination of substencils, so that higher order accuracy is recovered in smooth regions. In non-smooth regions, WENO reduces to making use of one of the substencils to obtain a properly winded, essentially non-oscillatory approximation. The notion of smoothness among each of the candidate substencils is accessed through a smoothness indicator that makes use of the first and second derivative information for the local approximations to identify non-smooth regions. Several other WENO schemes, which are more robust at shock capturing and generate less dissipation, have been proposed in the literature (see e.g., \cite{ABC,bor,hen,Serna}).  Traditional WENO schemes were originally based on polynomial interpolation.  However, algebraic polynomials are known to have limits in approximating data containing rapid gradients or fast oscillations due to their shift-and-scale invariant property. In order to address this problem, schemes were developed based on both trigonometric  \cite{ZQ1, ZQ2} and exponential  \cite{HKYY, HKYY-2} functions in the interpolation basis for ENO and WENO.

This paper proposes to use non-polynomial function approximations by formulating a WENO scheme in terms of RBFs and achieving superconvergence by tuning the available shape or tension parameter. RBFs are widely used as a basis function for multivariate scattered data approximation problems \cite{ RBF-B, olson2008gridless, RBF-Y}, and RBF approximation methods for solving partial differential equations have been developed in a variety of contexts \cite{RBF-SK, RBF-S}, including  WENO quadratures \cite{RBF-BH, GJ}.  A radial basis function $\phi:\RR^d \rightarrow \RR$ is defined in the sense that $\phi(x)=\phi(|x|)$, where $|\cdot|$ is the usual Euclidean norm.  Because of its definition, the power of the RBF approximation is in its meshless property, which is particularly beneficial in modeling scattered data. Moreover, the basis is flexible, as it can be tuned to incorporate local features of the data through the shape parameter.  For a given data set, an approximating function $Af(x)$ with an RBF $\phi$ can be represented as
\begin{align}\label{rbf-basic}
    Af(x)=\sum_{j=1}^M \alpha_j \phi (|x-\xi_j|)
\end{align}
where  $\{\xi_j : j=1,\cdots, M\}$ is a set of reference points and  $\alpha_j $ is a weight associated with $\phi(\cdot - \xi_j)$ for $j=1,\cdots, M$.
There are several ways to solve \eqref{rbf-basic} depending on the constraints. For example, assuming $f(\Omega)$ is given, then the constraints satisfy
\begin{align}\label{rbf-lin}
    Af(x^*)=f(x^*), \quad \forall x^*\in\Omega,
\end{align}
which may represent the solution of an interpolation problem (if $|\Omega| = M$) or an optimization problem ($|\Omega| \neq M$). In Table \ref{tab:rbf}, we provide commonly used global radial basis functions $\phi$ with a shape parameter $\lambda$. We note that the RBF approximation scheme  \eqref{rbf-basic} using a Gaussian function as a basis, i.e., $\phi(x)=e^{-(\lambda x)^2}$, has a conceptual resemblance with Gaussian Process (GP) modeling, which makes a probabilistic prediction instead of solving a linear system \eqref{rbf-lin}. In \cite{GP-R, GP-R2}, the authors use GP regression to solve hyperbolic conservation laws.

\begin{table}[t!]
	\centering
	\renewcommand{\arraystretch}{1.8}
	\begin{tabular}{l|c}
		\hline\hline
		RBFs &	$\phi(x)$   \\
		\hline
		Gaussian function & $\displaystyle e^{-(\lambda x)^2}$\\
		Multiquadric      & $\displaystyle \sqrt{1+(\lambda x)^2}$\\
		Inverse quadric   & $\displaystyle \frac{1}{1+(\lambda x)^2}$\\
		Thin-plate spline & $\displaystyle (\lambda x)^2 \log(\lambda x)$\\
		\hline\hline
	\end{tabular}
	\bigskip
	\caption{ Commonly used global radial basis functions $\phi$ with a shape parameter $\lambda$.}
	\label{tab:rbf}
\end{table}

The topic of superconvergence admits a vast array of literature. For example, several studies have explored the notion of superconvergence, with discontinuous Galerkin (dG) methods, on time evolution in ordinary differential equations \cite{DG1, DG2}, as well as hyperbolic and convection-diffusion PDEs \cite{DG3, DG4, DG5, DG6}. In dG schemes, superconvergent behavior can be incorporated into the basis by using information from the exact solution, provided one is available. Spectral methods for PDEs can also achieve high order accuracy, when the solutions are analytic, for continuous problems. For PDEs that admit discontinuous solutions, the Gibbs phenomenon is known to contaminate solutions, resulting in non-uniform convergence \cite{gegen0}.  However, in such instances, high order accuracy can be recovered through the use of certain post-processing techniques \cite{DG7, DG8, gegen1, gegen2}. In the literature, non-polynomial based numerical schemes have successfully demonstrated improvements in accuracy. We refer the interested readers to the papers \cite{GJ} and \cite{HKYY3}. The latter work develops a WENO scheme with a basis consisting of exponential polynomials, while the former proposes an ENO scheme that uses RBF interpolation. The primary objective of this paper is to devise high order FV schemes using compact, non-polynomial interpolation techniques, which achieve additional accuracy by exploiting the shape parameter available in the basis. As will be discussed in section \ref{sec: Superconvergent schemes}, the convergence order of the scheme is related to the level of accuracy of the approximation for shape parameter $\lambda^2$, which appears in the basis, e.g., $\phi(x)=e^{-(\lambda x)^2}$.

The organization of this paper is as follows. In section \ref{sec: RBF}, we begin with a general overview of RBF interpolation. Once we have introduced the interpolation problem, we demonstrate, in section \ref{sec: Superconvergent schemes}, how the shape parameter in the basis can be exploited to obtain superconvergent approximations. Then, in section \ref{sec: WENO schemes}, we briefly summarize the construction of FV WENO schemes for conservation laws, highlighting, in particular, subsection \ref{subsec: New smoothness indicators}, which defines the smoothness indicators employed by the proposed schemes. We then outline the key steps used in the implementation of the algorithms in section \ref{sec:newscheme}, which also contains details concerning extensions for two-dimensional problems. Experimental results, collected on a suite of test problems consisting of one-dimensional hyperbolic and weakly hyperbolic systems of conservation laws, as well as a two-dimensional nonlinear scalar test problem, are presented in section \ref{sec: Numerical results}. Finally, in section \ref{sec: Conclusion}, we briefly summarize the ideas presented in this work. 
	
\section{RBF Interpolation}\label{sec: RBF}

In this section, we provide a brief overview of interpolation with radial basis functions, which shall be useful for introducing our new WENO formulation. Suppose that a continuous function $f:\RR^d\rightarrow \RR$ is known only at a set of discrete points $X := \{x_1, . . . , x_N\}$ in $\Omega \subset \RR^d$. A function $\phi:\RR^d \rightarrow \RR$ is radial in the sense that $\phi(x)=\phi(|x|)$, where $|\cdot|$ is the usual Euclidean norm. RBF interpolation for $f$ on $X$ starts by choosing a basis function  $\phi$, and then defines an interpolant by
\begin{align}\label{rbf_interpol}
    Af_X(x):=\sum_{k=1}^{m} \beta_k p_k(x)+\sum_{j=1}^{N}\alpha_j\phi(x-x_j),
\end{align}
where $\{p_1, \ldots, p_m \}$ is a basis for $\Pi_m $ and the coefficients $\alpha_j$ and $\beta_i$ are chosen to satisfy the linear system
\begin{align}\label{rbf_lin}
\begin{cases}
\begin{split}
    Af_X(x_i) = \sum_{k=1}^{m} \beta_k p_k(x_i)+\sum_{j=1}^{N}\alpha_j\phi(x_i-x_j) &= f(x_i), \quad i=1,\cdots,N,\\
    \sum_{j=1}^{N} \alpha_j p_k(x_j) &= 0, \quad \qquad k=1,\cdots,m.
\end{split}
\end{cases}
\end{align}
Here $\Pi_m$ denotes the space generated by all algebraic polynomials of degree less than $m$ on $\RR^d$.
For a wide choice of functions $\phi$ and polynomial orders $m$, the existence and uniqueness of the solution of the linear system \eqref{rbf_lin} is ensured when $\phi$ is a conditionally positive definite function.   
\begin{definition}
   Let  $\phi:\RR^d \rightarrow \RR$ be a continuous function. We say that $\phi$ is conditionally positive definite of order  $m\in\NN$ if for every finite set of pairwise distinct points $X = \{x_1, \ldots x_N\} \subset \RR^d$ and $\alpha = (\alpha_1, \ldots ,\alpha_N) \in \RR^N \backslash \{0\} $ satisfying
   $\sum_{j=1}^{N} \alpha_j p(x_j) = 0$ for $\forall p\in\Pi_m$,
   the quadratic form $$\sum_{i=1}^N \sum_{j=1}^N \alpha_i\alpha_j \phi(x_i-x_j)$$
   is positive definite.
\end{definition}
This leads to the linear system \eqref{rbf_lin} for  $\alpha = (\alpha_1, \ldots ,\alpha_N) $ and  $\beta = (\beta_1, \ldots ,\beta_m) $, which is given in block-matrix form as
\begin{gather}\label{rbf_mat}
 \begin{bmatrix} {\Phi} & P^T \\ P & O \end{bmatrix}
 \begin{bmatrix} \alpha^T \\ \beta^T \end{bmatrix}
 =
 \begin{bmatrix} \mathbf{f}^T \\ O \end{bmatrix}
\end{gather}
where ${\Phi}=\{\phi(x_i-x_j):i,j\}$, $P=\{p_k(x_j):k,j\}$ and $\mathbf{f}=\{f(x_i):i\}$. If we assume $m=0$ in \eqref{rbf_interpol}, then from the equations \eqref{rbf_interpol} and \eqref{rbf_mat}, the interpolant can be represented as
\begin{align}\label{rbf_kernel}
    Af_X(x)=\sum_{j=1}^{N}\alpha_j\phi(x-x_j)= {\upphi}{\Phi}^{-1}\mathbf{f}^T,
\end{align}
where ${\upphi} = \{\phi(x-x_j):j\} $. We note that the product ${\upphi}{\Phi^{-1}}$, which appears as part of the representation \eqref{rbf_kernel}, is completely independent of the function values $\mathbf{f}$. Next, in section \ref{sec: Superconvergent schemes}, expressions for the optimal shape parameters, i.e., $\lambda^2$ in the non-polynomial basis, are derived, which results in superconvergent approximations.

	\section{Supercovergent RBF schemes}
\label{sec: Superconvergent schemes}

In this section, we derive the expressions for optimal shape, or tension, parameters made available through RBF interpolation. We perform our analysis in section \ref{subsec: Parameter} based on two reconstruction methods: a direct approach using integrals and a second approach that utilizes a primitive function. Here, we define the shape parameter $\lambda^2$, for a given RBF $\phi$, so that it maximizes the convergence order of the approximation. Once the techniques have been demonstrated, we generalize these results in sections \ref{subsec:N_odd} and \ref{subsec:N_even}, with the key points summarized in Theorems \ref{thm:odd} and \ref{thm:even}, respectively.

\subsection{Optimal Shape Parameters for RBFs}\label{subsec: Parameter}

In \cite{HKYY}, the authors introduced a WENO scheme based on the space of exponential polynomials. Later, in the work \cite{HKYY-2}, they improved the order of accuracy of their schemes by exploiting the control parameter $\lambda \in \RR$ or $i\RR$ for exponential basis functions of the form $e^{\lambda x}$. We adopt a similar strategy in this work with the difference being the choice of the basis. Here, the basis functions consist of infinitely smooth RBFs $\phi$ (see Table \ref{tab:rbf}) rather than exponential polynomials. We present the analysis using the Gaussian function
\begin{align}\label{rbf:gauss}
    \phi(x)=e^{-\lambda^2 x^2}, \quad \lambda\in \RR \text{ or } i\RR,
\end{align}
in the case of $N = 2$, using two different approaches, but these techniques can be easily extended to other RBFs (see Table \ref{tab:rbf} for other options). For each approach, we provide analytical expressions for the optimal shape parameters, which allow the schemes to achieve superconvergence on a fixed reconstruction stencil. 



\subsubsection{Direct Computation with Integrals}
\label{subsubsec: Direct with integrals}


Recall that the RBF approximation is given by
\begin{align}\label{rbf_app2}
Au(x):= \sum_{k=1}^{2} \alpha_k  \phi(x-x_{j+k-1})= \alpha_1  \phi(x-x_{j}) +\alpha_2  \phi(x-x_{j+1}).
\end{align}
The goal is to form a high order approximation using the form \eqref{rbf_app2}, which preserves each of the cell averages. In other words, the approximation should satisfy the integral constraint
\begin{align}\label{rbf_app2_cond}
    \frac{1}{\Delta x}\int_{I_i} Au(\xi) d\xi =\bar{u}_i, \quad i=j,{\ }j+1.
\end{align}
Using the RBF $\phi(x)$ defined in \eqref{rbf:gauss}, it follows that the integrals of $\phi$ can be evaluated analytically, which resulting in a solvable linear system. Once the solution is determined, the final approximation at the cell boundary $x_{j+\frac{1}{2}}$ is given by
\begin{align*}
    Au(x_{j+\frac{1}{2}})=\frac{ 2\lambda \Delta x\exp(-\frac{\lambda^2\Delta x^2}{4}) }
    { \sqrt{\pi} \big( \erf(\frac{\lambda \Delta x}{2}) + \erf(\frac{3\lambda\Delta x}{2})  \big)  } (\bar{u}_j+\bar{u}_{j+1}).
\end{align*}
To determine the shape parameter, we Taylor expand the right-hand side of the previous equation, which yields
\begin{align*}
    Au(x_{j+\frac{1}{2}}) &\approx \left(\frac{1}{2} +\frac{\lambda^2\Delta x^2 }{6} - \frac{7\lambda^4\Delta x^4}{90}\right)(\bar{u}_j+\bar{u}_{j+1}) \\
&= u(x_{j+\frac{1}{2}}) + {\Delta x^2}\Bigg(\frac{1}{3} {\lambda^2}u(x_{j+\frac{1}{2}}) + \frac1 6 u''(x_{j+\frac{1}{2}}) \Bigg) + \mathcal{O}(\Delta x^4).
\end{align*}
Hence, we have a second order approximation for the cell boundary point. Notice that we can obtain a fourth order approximation if we choose the shape parameter with
\begin{align*}
   {\lambda^2}=-\frac{u''(x_{j+\frac{1}{2}})}{2u(x_{j+\frac{1}{2}})}+\mathcal{O}(\Delta x^2).
\end{align*}

In practice, this parameter can be formed entirely from the available cell average data. Although there may be many ways to approximate these terms, our implementation considers finite difference approximations to the derivatives appearing in the expression for the optimal shape parameter. We provide additional details in Appendix \ref{app:lambda}. Note that the cell average data can be shifted by some positive constant at the beginning and end of the reconstruction steps to prevent division by zero, if needed. In the next section, we perform a similar reconstruction using a primitive function.




\subsubsection{Construction with a Primitive Function}
\label{subsubsec: A primitive function}

To reconstruct the approximation at the cell boundary $x=x_{j+\frac{1}{2}}$ from the cell averages $\{\bar{u}_j,\bar{u}_{j+1}\}$, we define a primitive function $$U(x)=\int_{x_{j-\frac{1}{2}}}^x u(\xi) d\xi$$ which can be explicitly written in terms of the available cell averages as 
$$U(x_{i-\frac{1}{2}})= \Delta x \sum_{\ell =j}^{i-1} \bar{u}_{\ell} \quad \text{ for } i=j,\cdots,j+2.$$
Using the available interpolatory data for the primitive function at the cell interfaces, we seek an RBF representation 
\begin{align}\label{rbf_app1}
AU(x):= \sum_{k=0}^{2} \alpha_k \phi \left(x-x_{j+k-\frac{1}{2}}\right),
\end{align}
which satisfies
\begin{align}\label{rbf_app1_cond}
    AU(x_{i-\frac{1}{2}})=U(x_{i-\frac{1}{2}}), \quad i=j,j+1,j+2.
\end{align}
Then the final approximation to $u(x)$ at $x=x_{j+\frac{1}{2}}$ is obtained by differentiating the RBF representation
$$AU'(x):= \sum_{k=0}^{2} \alpha_k \phi'(x-x_{j+k-\frac{1}{2}}),$$ which approximates $U'(x)=u(x)$, i.e., 
\begin{align}\label{rbf_app1_fin0}
    AU'(x_{j+\frac{1}{2}}) \approx U'(x_{j+\frac{1}{2}})=u(x_{j+\frac{1}{2}}).
\end{align}
The approximation with a Gaussian RBF $\phi(x)$ is found to be
\begin{align}\label{eq:rbf_loc}
    AU'(x_{j+\frac{1}{2}})=2{\lambda^2}\Delta x^2\frac{ e^{3{\lambda^2}\Delta x^2} }{ e^{4{\lambda^2}\Delta x^2}-1}(\bar{u}_j+\bar{u}_{j+1}),
\end{align}
with the basis coefficients calculated from \eqref{rbf_app1_cond} and \eqref{rbf_app1_fin0}.
Applying a Taylor expansion to the right-hand side of this equation gives 
\begin{align*}
    AU'(x_{j+\frac{1}{2}}) &\approx \left(\frac{1}{2} +\frac{1}{2}\lambda^2\Delta x^2 - \frac{1}{12}{\lambda^4\Delta x^4}\right)(\bar{u}_j+\bar{u}_{j+1}) \\
&=u(x_{j+\frac{1}{2}}) + \Delta x^2 \Bigg({\lambda^2} u(x_{j+\frac{1}{2}})+\frac1 6 u''(x_{j+\frac{1}{2}})\Bigg) + \mathcal{O}(\Delta x^4).
\end{align*}
Hence, we have the second order approximation for the reconstructed cell boundary point. Alternatively, we can obtain the fourth order approximation if we choose the shape parameter with
\begin{align}\label{lam:local}
   {\lambda^2}=-\frac{u''(x_{j+\frac{1}{2}})}{6u(x_{j+\frac{1}{2}})}+\mathcal{O}(\Delta x^2),
\end{align}
which can be computed using the techniques explained at the end of section \ref{subsubsec: Direct with integrals} and further described in Appendix \ref{app:lambda}. Using this same machinery, we can generalize these approximations by considering the parity of $N$, which results in two theorems, which seek to address the overall convergence order.

\subsection{Three-Point RBF Schemes}
\label{subsec:N_odd}

Following the analysis of the previous section, we construct RBF interpolation schemes with optimal parameters when the number of the stencil points $N$ is odd. Through a fairly direct construction, one obtains the following theorem:

\begin{theorem}\label{thm:odd}
	Let $u$ be a smooth function on $\Omega$ and $\phi$ be a smooth radial basis function. Given a set of reference points $\{x_k\in\Omega:k=1,\cdots,N\}$ for an odd integer $N$, there exists a set of coefficients \hbox{$\{\alpha_k :k=1,\cdots,N\}$} of the approximation
	\begin{align*}
	Au(x) = \sum_{k=1}^{N} \alpha_k \phi(x-x_k)
	\end{align*}
	to the function $u(x)$, which are constructed from cell averaged data $\{\bar{u}(x_j) : x_j \in \Omega\}$. Furthermore, the approximation can be made $(N+1)$st order accurate, i.e.,
	\begin{align*}
	Au(x)=u(x)+ \mathcal{O}(\Delta x^{N+1}),
	\end{align*}
	for $x\in\Omega$.
\end{theorem}

We show an example for the case of $N=3$ on behalf of the proof of Theorem \ref{thm:odd}. Suppose that we construct the approximation of a smooth function $u$ at $x=x_{j+\frac12}$ using the cell averages \mbox{ $\{\bar{u}_{j}:j=-1,0,1\}$}. Using a primitive function 
$$U(x)=\int_{x_{j-\frac{1}{2}}}^x u(\xi) d\xi, $$
the approximation using a radial basis function $\phi(x)=e^{-\lambda^2 x^2}$ is defined by
\begin{align*}
AU(x):= \sum_{k=-1}^{2} \alpha_k \phi(x-x_{j+k-\frac12}),
\end{align*}
so that $AU'(x)\approx u(x)$.
Repeating the steps outlined in subsection \ref{subsubsec: A primitive function}, and Taylor expanding the RBF approximation, we find that
\begin{align*}
    AU'(x_{j+\frac12}) 
    &= u(x_{j+\frac12}) + {\Delta x^3}\left( {\lambda^2}u'(x_{j+\frac12}) + \frac{1}{12}u'''(x_{j+\frac12}) \right) + \mathcal{O}(\Delta x^4).
\end{align*}
Therefore the third order scheme can be improved to fourth order accuracy through the choice
\begin{align}\label{eq:lam4rbf3}
    \lambda^2= -\frac {u'''(x_{j+\frac12})}{12u'(x_{j+\frac12})} +  \mathcal{O}(\Delta x),
\end{align}
which can be computed using the same techniques discussed in Appendix \ref{app:lambda}.

\begin{remark}\label{rmk:odd effective stencil}
    Following the brief discussion in Appendix \ref{app:lambda}, the approximation \eqref{eq:lam4rbf3} requires one \mbox{additional} cell average value that lies outside of the three-point global reconstruction stencil for $u_{j+\frac{1}{2}}^{-}$. Since the derivatives in this approximation do not need to be winded, this cell average value can be selected to lie in the union of global reconstruction stencils used to form $u_{j+\frac{1}{2}}^{-}$ and $u_{j+\frac{1}{2}}^{+}$, i.e., $\{ \bar{u}_{j-1}, \bar{u}_{j}, \bar{u}_{j+1}, \bar{u}_{j+2} \}$. This is identical to the stencil used by the third order WENO-JS scheme when forming the flux at $x_{j+1/2}$.
\end{remark}

\subsection{Four-Point RBF Schemes}
\label{subsec:N_even}

Here we provide the following theorem for the case when $N$ is even to complete the analysis of proposed scheme. Again, through a fairly straightforward construction, one obtains the general theorem, which is as follows:

\begin{theorem}\label{thm:even}
	Let $u$ be a smooth function on $\Omega$ and $\phi$ is a smooth radial basis function. Given a set of reference points $\{x_k\in\Omega:k=1,\cdots,N\}$ for an even integer $N$, there exists a set of coefficients $\{\alpha_k :k=1,\cdots,N\}$ of the approximation
	\begin{align*}
	Au(x) = \sum_{k=1}^{N} \alpha_k \phi(x-x_k)
	\end{align*}
	to the function $u(x)$, which are constructed from cell averaged data $\{\bar{u}(x_j) : x_j \in \Omega\}$. Moreover, the resulting approximation is $(N+p)$th order accurate, with $0\leq p\leq 2$, i.e.,
	\begin{align*}
	Au(x)=u(x)+ \mathcal{O}(\Delta x^{N+p}),
	\end{align*}
	for $x\in\Omega$.
\end{theorem}

Here, we provide the analysis for the case $N=4$. Using the cell averages $\{\bar{u}_{j}:j=-1,\cdots,2\}$, the approximation is defined by
\begin{align}\label{rbf_app1_4}
AU(x):= \sum_{k=-1}^{3} \alpha_k \phi(x-x_{j+k-\frac{1}{2}})
\end{align}
and at the cell boundary $x=x_{j+\frac{1}{2}}$, we have that
\begin{align}\label{rbf_app1_4_cond}
\begin{split}
    AU(x_{i-\frac{1}{2}})&=U(x_{i-\frac{1}{2}})\\ 
    &= \Delta x \sum_{\ell =j}^{i-1} \bar{u}_{\ell} 
    , \quad i=j-1,\cdots,j+3.
\end{split}
\end{align}
Proceeding as before, we can obtain the solution 
\begin{align}\label{eq:rbf_global}
     AU'(x_{j+\frac{1}{2}}) &\approx \sum_{k=-1}^{2} C_k \bar{u}_{j+k} 
\end{align}
with coefficients $C_k$ computed from equations \eqref{rbf_app1_4} and \eqref{rbf_app1_4_cond}. As before, these
can be Taylor expanded about the cell boundary $x=x_{j+\frac{1}{2}}$ and we find that
\begin{align*}
    AU'(x_{j+\frac{1}{2}})
    &= u(x_{j+\frac{1}{2}}) + {\Delta x^4}\left(\frac{1}{30}u^{(4)}(x_{j+\frac{1}{2}}) +  \frac23{\lambda^2}u''(x_{j+\frac{1}{2}}) + 2 u(x_{j+\frac{1}{2}}){\lambda^4} \right) + \mathcal{O}(\Delta x^6).
\end{align*}
Next, we choose $\lambda^2$ to remove the  $\mathcal{O}(\Delta x^4)$ term, i.e.,
\begin{equation}\label{eq:lam_choice}
    \lambda^2=\frac { -\frac{1}{3}u''(x_{j+\frac{1}{2}}) \pm \sqrt{\frac1 9 u''(x_{j+\frac{1}{2}})^2 - \frac{1}{15} u(x_{j+\frac{1}{2}})u^{(4)}(x_{j+\frac{1}{2}})}  } {2u(x_{j+\frac{1}{2}})} +  \mathcal{O}(\Delta x^p),
\end{equation}
which can be computed using the same techniques discussed in Appendix \ref{app:lambda}. Therefore we can obtain a superconvergent scheme, which is $(4+p)$th order accurate with $p\leq 2$. We now proceed to the discussion of WENO schemes in section \ref{sec: WENO schemes}, making use of the superconvergent RBF approximations.

\begin{remark}\label{rmk:even effective stencil}
    In contrast to Remark \ref{rmk:odd effective stencil}, evaluating the approximation \eqref{eq:lam_choice}, with $p=2$, requires one \mbox{additional} cell average value that lies outside of the union of global reconstruction stencils, in the fifth order WENO-JS scheme, to compute $u_{j+\frac{1}{2}}^{-}$ and $u_{j+\frac{1}{2}}^{+}$. On the other hand, when $p=1$, the ``effective stencil" is no larger than the global stencil used by the fifth order WENO-JS scheme in the construction of the flux at $x = x_{j+\frac{1}{2}}$.
\end{remark} 

\section{WENO Schemes}\label{sec: WENO schemes}

This section describes the general formulation of FV WENO schemes used to solve conservation laws. First, we provide a brief summary of the FV discretization in section \ref{subsec: FV schemes}. Then, in section \ref{subsec: New smoothness indicators}, we introduce new smoothness indicators, which were motivated by numerical experimentation, and discuss the mapping for the WENO weights.

\subsection{Formulation of a Finite Volume Scheme}\label{subsec: FV schemes}

Consider one-dimensional variant of \eqref{eq:claw0}, which takes the form 
\begin{align}
\label{eq:claw1}
\begin{split}
     & u_t + f(u)_x  = 0,\quad   \ x\in \RR,\\
     &     u(x,0)    = u_0(x).
\end{split}
\end{align}
To develop a FV discretization, we let the computational domain be partitioned into uniform cells, so that the $j$th cell is given by $I_j=[x_{j-1/2}, x_{j+1/2}]$. Further, since the cells are uniformly spaced, each has the size $\Delta x = x_{j+1/2} - x_{j-1/2}$. FV schemes require cell averages of the solution $u$
\begin{align}\label{eq: cell averages}
\bar{u}_j(t) :=\frac{1}{\Delta x}\int_{I_j} u(x,t) dx,
\end{align}
where the integrals in \eqref{eq: cell averages} can be numerically approximated through quadrature of suitable accuracy. By integrating the \eqref{eq:claw1} over the control volume $I_j$, one obtains a collection of evolution equations defined in each of the control volumes, i.e.,
\begin{equation}\label{SEMI-DIS}
	\frac{d}{d t} \bar{u}_j(t) = - \frac{1}{\Delta x} \left( f(u(x_{j+1/2})) - f(u(x_{j-1/2})) \right).
\end{equation}
Defining the numerical flux $\hat{f}_{j\pm 1/2}$ by
\begin{align}\label{num_flux}
    \hat{f}_{j\pm 1/2} = h(u^-_{j\pm 1/2}, u^+_{j\pm 1/2}),
\end{align}
the equation \eqref{SEMI-DIS} is approximated as 
\begin{equation}\label{eq_app}
\frac{d}{d t} \bar{u}_j(t) = - \frac{1}{\Delta x} \left( \hat{f}_{j + 1/2} - \hat{f}_{j - 1/2} \right).
\end{equation}
The monotone flux $h$ satisfies several properties, namely, it is Lipschitz continuous in both arguments and should be consistent with the physical flux $f$, i.e., $h(u,u)=f(u)$. Moreover, the flux should be non-decreasing (non-increasing) with respect to the first (second) argument. In this paper, we employ the HLLC and Lax-Friedrichs fluxes for solving the hyperbolic problems and a Godunov flux for the weakly hyperbolic system. The definitions of these numerical fluxes can be found in section \ref{sec: Numerical results}. Next, we focus on the WENO component of the proposed schemes, which seeks to develop high order reconstructions for the cell average data supplied to the numerical flux functions, i.e., $u^-_{j\pm 1/2}$ and $ u^+_{j\pm 1/2}$.  

\begin{figure}
    \centering
    \includegraphics{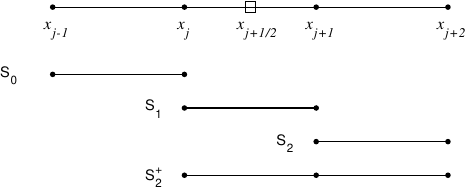} 
    \caption{Three two-point substencils $S_k$, $k=0,1,2$ and a three-point substencil $S^+_{2}$.}
    \label{fig:stencil}
\end{figure}


\subsection{Construction of New Smoothness Indicators for WENO-RBF Schemes}
\label{subsec: New smoothness indicators}

It is well known that the smoothness indicator plays a pivotal role in the WENO reconstruction procedure. 
We are proposing a three-point WENO scheme and a four-point WENO scheme based on RBFs with optimal parameters, so it will be helpful to introduce some references to simplify names for the schemes. From now on, we shall refer to the three-point and four-point WENO schemes, which are based on RBFs, as WENO-RBF3 and WENO-RBF4, respectively. Furthermore, the WENO-RBF3 scheme is constructed using two two-point substencils ($S_0$, $S_1$ in Figure \ref{fig:stencil}), while the WENO-RBF4 scheme is using three two-point substencils ($S_0, S_1, S_2$ in Figure \ref{fig:stencil}) so that only the first-order difference operator can be used to estimate the smoothness of numerical fluxes. In other words, the small size of the substencils, used by the RBF methods, introduces limitations which make capturing highly oscillatory structures and rapid gradients challenging, especially on coarse grids. Through experimentation, we also found that switching to the smoothness indicators used in classical WENO (i.e., WENO-JS) resulted in minimal improvement.

In an effort to amend this issue, we employ difference operators which have so-called \textit{exponentially vanishing moments}. The basic idea is that difference operators are constructed so that their applications to smooth functions result in more rapid convergence to zero than classical undivided differences. Such operators should be more effective at detecting smoothness, or lack thereof, in the data, even on a small collection of points. Before proceeding, we first introduce some ideas from a recent work \cite{HKYY3}, in which smoothness indicators were constructed using exponential polynomials.

Exponential polynomials can be characterized as the kernel of an operator  $p_m(D)$ of the form
\begin{equation} \label{PD-10}
    p_m(D) f := \prod_{i=1}^m (D-\gamma_i I) f,
\end{equation}
where $D$ is a continuous  differential operator, $I$ is an identity operator and $\gamma_i \in \CC$, which is associated with the target exponential polynomials.
Then, the operator $p_m(D) f$ is constructed as a differential operator which annihilates exponential polynomials. For our purposes, it suffices to consider the case $m=1$ in equation \eqref{PD-10}. To illustrate, consider the exponential function $f(x)=e^{-\gamma x}$ for some constant $\gamma$.
Then, the first order operator $p_1(D)$ with the form  $ p_1(D) f  = (D +\gamma  I)f$ annihilates the exponential function, i.e., 
\begin{equation*}
    f(x) = e^{-\gamma x} \implies {p}_1(D) f= 0.
\end{equation*}
For the discretized version, we assume the function values $\{f(x_{i-1}), f(x_i)\}$ on the stencil $\{x_{i-1}, x_i\}$ are available. Then, the discrete version of  $p_1(D) f$ is defined as
$$
{\bf p}_1(D)_i f :=(D+\gamma I)_i f:= e^{\gamma\Delta x} f(x_{i})  - f(x_{i-1}),
$$
which also annihilates the exponential function. Similar constructions can be achieved for $m\geq 2$, leading to $m$th order operators, which annihilate $m$ exponential polynomials. For additional information, we refer interested readers to \cite{HKYY3}.



We are now ready to propose new smoothness indicators based on the new undivided difference with exponentially vanishing moments. Specifically, we define $\beta_k$, $k=0,1,2$, by
\begin{align}\label{BETA-Y}
\begin{split}
\beta_0 &:=| D_{i-1}^1 f |^2 + | {\bf p}_1(D)_{i-1} f |^2
        = | f_{i} -f_{i-1} |^2 + |e^{\gamma_i\Delta x} f_{i} -f_{i-1} |^2, \\
\beta_1 &:=| D_i^1 f |^2 + | {\bf p}_1(D)_{i} f |^2
          = | f_{i+1} - f_i|^2 + |e^{\gamma_{i+1}\Delta x} f_{i+1} - f_i|^2, \\
\beta_2 &:=\frac1 2 \big( \beta_1 + (| D_{i+1}^1 f |^2 + | {\bf p}_1(D)_{i+1} f |^2) \big) \\
     & = \frac1 2 \big(  | f_{i+1} - f_i|^2 + |e^{\gamma_{i+1}\Delta x} f_{i+1} - f_i|^2 + | f_{i+2} - f_{i+1}|^2 + |e^{\gamma_{i+2}\Delta x} f_{i+2} - f_{i+1}|^2  \big).
\end{split}
\end{align}
In smooth regions, $\beta_k$ should be small, so we select the parameter $\gamma_i$ in a way that depends on the given stencil data:
\begin{equation}\label{GAM-20}
\gamma_{i+\nu}\Delta x =- \frac{f_{i+\nu} -f_{i+\nu-1}} {f_i +\delta},  \qquad
{\rm sign} (\delta) := {\rm sign} (f_i), \quad \nu \in \{0,1,2\}.
\end{equation}
where $\delta=\delta(\Delta x)$ is introduced to prevent the denominator from becoming zero.
Observe that $\beta_2$ is defined using three points ($S_2^+$ in Figure \ref{fig:stencil}) instead of two points, which incorporates a bias in the indicators.

Next, using the local smoothness indicators, we map the linear weights $d_k$ to the nonlinear weights $\alpha_k$ for $k=0,1,2$ via
\begin{align}\label{ALP}
\alpha_k = d_k \Bigg( 1 + \frac{\tau}{\beta_k+ \varepsilon } + \left(\frac{\beta_k}{\tau+ \varepsilon }\right)^2 \Bigg),
\quad \varepsilon: = \varepsilon(\Delta x),
\end{align}
where $\epsilon > 0 $ is used to prevent the denominator from becoming zero.
Here $\tau$ measures the global smoothness and is defined by $ \tau = |\beta_2 - \beta_0|$ in WENO-RBF3 and $ \tau = |\beta_1 - \beta_0|$ in WENO-RBF4.
The linear weights $\{d_k:k=0,1,2\}$ are chosen so that the linear combination of fourth order local approximations, on each of the two-point stencils, is consistent with the sixth order approximation obtained on the four-point stencil. The nonlinear weights for WENO-RBF3 are then scaled to form a partition of unity 
\begin{align}\label{OMEGA3}
\omega_k = \frac{\alpha_k}{\sum_{\ell=0}^1 \alpha_{\ell} }, \quad
k = 0,1
\end{align}
and the nonlinear weights for WENO-RBF4 are defined by
\begin{align}\label{OMEGA}
\omega_k = \frac{\alpha_k}{\sum_{\ell=0}^2 \alpha_{\ell} }, \quad
k = 0,1,2.
\end{align}
The value $u_{j+1/2}^{-}$ is approximated by a convex combination of local approximations over each of the substencils $S_k$ using the nonlinear weights $\omega_k$ so that
$$u_{j+1/2}^{-} := \sum _{k=0}^{1} \omega_k u^{(k)}_{j+1/2}, $$
and
$$u_{j+1/2}^{-} := \sum _{k=0}^{2} \omega_k u^{(k)}_{j+1/2}, $$
in the WENO-RBF3 scheme and WENO-RBF4 scheme, respectively. The analogous construction for $u_{j+1/2}^{+}$ follows by symmetry. In the next section, we present our full FV WENO-RBF algorithm and some details concerning the implementation.
	
\section{Implementing Superconvergent WENO-RBF Schemes}\label{sec:newscheme}

Now that we have introduced the key components of the proposed schemes (see sections \ref{sec: Superconvergent schemes} and \ref{sec: WENO schemes}), we can describe the implementation. First, we begin with some details regarding hybrid WENO schemes in section \ref{subsec: Hybrid scheme}, before summarizing the key steps of the FV WENO-RBF algorithm in section \ref{subsec: FV RBF-WENO Algorithm}. Then, in section \ref{subsec: twoD-scheme}, we briefly discuss the implementation of the proposed schemes to multi-dimensional problems, focusing, in particular on the two-dimensional scalar case. Stencil coefficients for the RBF (Gaussian) methods developed in this work can be found in Appendix \ref{app:RBF_coeffs}.

\subsection{Comments on the Hybrid Implementation}
\label{subsec: Hybrid scheme}

The implementation of the WENO-RBF4 method used in this work employs a hybrid strategy, which aims to alleviate the computational cost associated with WENO methods due to the additional cost from smoothness indicators and mappings for the nonlinear weights. The basic idea of a hybrid approach is to use reconstructions on a fixed set of cell average data in regions where the data is smooth, while non-smooth regions are appropriately handled with a WENO scheme (see e.g., section \ref{sec: WENO schemes}). The adaptive selection of a reconstruction method relies on the use of certain smoothness criterion,  which are, ideally, inexpensive to evaluate. This criterion can, for example, be evaluated with smoothness indicators, such as those in the classical WENO approach \cite{Jiang1}, as well as divided or undivided differences. A tolerance (or threshold) then selects the reconstruction method according to the smoothness of the given data. Our selection process consists of the following steps:
\begin{enumerate}
    \item Using finite differences, first compute the relative smoothness $r(x_i)$, which we define as 
    \begin{equation}
        \label{eq:relative smoothness}
        r(x_i) = \frac{2 \Big( \Big\lvert \delta [\bar{u}](x_{i}) \Big\rvert + \Big\lvert \delta^{2} [\bar{u}](x_{i}) \Big\rvert \Big) }{ \Big\lvert \Delta^{-} [\bar{u}](x_{i-1}) \Big\rvert + \Big\lvert \delta^{2} [\bar{u}](x_{i-1}) \Big\rvert + \Big\lvert \Delta^{+} [\bar{u}](x_{i+1}) \Big\rvert + \Big\lvert \delta^{2} [\bar{u}](x_{i+1}) \Big\rvert }, \quad i = 1, \cdots, N,
    \end{equation}
    where $\delta$ and $\delta^2$ are central difference operators for first and second derivatives, and $\Delta^{-}$ and $\Delta^{+}$ are the backward and forward difference operators for the first derivative. These difference approximations are all computed to second order accuracy. Note that along the boundaries, data from an extension is required and can be constructed using extrapolation of sufficient accuracy. 
    \item Once step 1 is complete, we find the minimum and maximum values of the relative smoothness $r_{min}$ and $r_{max}$ and then compute the tolerance  
    \begin{equation}
        \label{eq:tolerance for the hybrid scheme}
        r_{tol} = \min \left( \theta, \kappa \frac{ r_{min} + \epsilon }{ r_{max} + \epsilon }\right).
    \end{equation}
    In the numerical experiments which use the hybrid approach, we take $\theta = 1.5$, $\kappa = 5$, and $\epsilon = 1 \times 10^{-10}$. 
    \item Next, we map $r(x_{i}) \mapsto \{0,1\}$ using the previously computed tolerance $r_{tol}$. We identify cells, which are to use WENO reconstructions, as those for which $r(x_i) \geq r_{tol}$, using fixed stencil reconstructions for those that remain. To account for shortcomings in the definition of the relative smoothness \eqref{eq:relative smoothness}, we flag a buffer zone of 4 cells in each direction around any cell marked for WENO reconstruction.    
\end{enumerate}

The parameter choices used for $\theta$ and $\kappa$ are based on numerical experimentation and are, by no means, exhaustive. While it is entirely possible that better choices for these parameters exist, it is beyond the scope of this work. Experimentation with different parameters mostly resulted in a more conservative hybrid algorithm, where the WENO reconstruction process was being applied in larger regions of the domain, even where the solution is smooth. While these selections were not particularly detrimental to the shock-capturing abilities of the method, the resulting imbalance does increase the time-to-solution, as the WENO reconstruction is more expensive.

    
    
    
    
    
    
    
    
    
    
    
    

\subsection{The FV WENO-RBF Component}
\label{subsec: FV RBF-WENO Algorithm}


Once a method for each cell has been selected, the algorithm then applies the corresponding reconstruction technique. The basic steps in the algorithm for the WENO-RBF3 and WENO-RBF4 methods are exactly the same, so we will only present a summary for WENO-RBF4. As discussed in the previous section, cells deemed smooth use a four-point fixed stencil RBF interpolant; regions characterized as ``non-smooth" apply the WENO-RBF algorithm, which proceeds as follows:

\begin{enumerate}
\setcounter{enumi}{0}
    \item Using \eqref{eq:rbf_loc}, form the local approximations $u^{(k)}_{j+1/2}$ on each of the two-point substencils $S_k:=\{x_{j+k-1},x_{j+k}\}$, with $k = 0,1,2$: 
   \begin{align*}
          u^{(k)}_{j+1/2} := \sum_{\ell =0}^{1} c^k_{\ell} \bar{u}_{j+k-1+\ell}.
   \end{align*}
   The coefficients $c_{\ell}^k$, where $\ell=0,1$ and  $k=0,1,2,$ are defined with the local shape parameter for the RBF $\phi(x)=\exp(-\lambda_L^2 x^2)$ obtained from \eqref{lam:local}:
    \begin{align*} 
        {\lambda_L^2}\approx -\frac{u''(x_{j+\frac{1}{2}})}{6u(x_{j+\frac{1}{2}})}.
   \end{align*}
   Here the function values are replaced with cell averages and the derivatives are obtained with finite differences.
   \item Using \eqref{eq:rbf_global}, we can form the approximation $u^{S_4}_{j+1/2}$ on the big stencil $S_4:=\{x_{j-1}, \cdots, x_{j+2}\}$ as
   \begin{align*}
          u^{S_4}_{j+1/2} = \sum_{\ell=-1}^{2} C_{\ell} \bar{u}_{j+\ell}, 
   \end{align*}
   where $C_{\ell}$, and $\ell=-1, \cdots, 2$ is a coefficient which is dependent on the global shape parameter $\lambda_G$. For the RBF $\phi(x)=\exp(-\lambda_G^2 x^2)$, this is reflected in \eqref{eq:lam_choice}:
    \begin{align*} 
        {\lambda_G^2}\approx \frac { -\frac{1}{3}u''(x_{j+\frac{1}{2}}) \pm \sqrt{\frac1 9 u''(x_{j+\frac{1}{2}})^2 - \frac{1}{15} u(x_{j+\frac{1}{2}})u^{(4)}(x_{j+\frac{1}{2}})}  } {2u(x_{j+\frac{1}{2}})},
   \end{align*}
   As in step 1, the function values are replaced with cell averages and the derivatives are obtained with finite differences, as described at the end of section \ref{subsubsec: Direct with integrals}.
   \item Compute the linear WENO weights $\{d_k:k=0,1,2\}$ which satisfy 
   $$\sum _{k=0}^2 d_k u^{(k)}_{j+1/2} = u^{S_4}_{j+1/2} + O(\Delta x^p),$$
   using an appropriate high order $p$ from step 1 and step 2. This results in weights of the form 
   $$d_0 = \frac{C_{-1}}{c_{0}^0}, \quad d_2 = \frac{C_2}{c_1^2}, \quad d_1 = 1-d_0 - d_2.$$ This reflects the partition of unity for the linear weights.
    \item Using the cell averages, compute the parameters \eqref{GAM-20} and the smoothness indicators \eqref{BETA-Y}.
    \item Map the linear weights $\{d_k:k=0,1,2\}$ from step 3 to nonlinear weights $\{\omega_k:k=0,1,2\}$ using equations \eqref{ALP} and \eqref{OMEGA}.
    \item Obtain the reconstructed value at the cell interface using the nonlinear weights: 
    $$u_{j+1/2}^{-} = \sum _{k=0}^{2} \omega_k u^{(k)}_{j+1/2} .$$
\end{enumerate}
The procedure for determining $u_{j+1/2}^{+}$ follows, analogously, by reflecting the cell average stencil data. Once the reconstructions for $u_{j+1/2}^{\pm}$ are completed, we simply apply the numerical flux function \eqref{num_flux}, which yields $\hat{f}_{j+ 1/2}$. In the case of hyperbolic systems, these reconstructions are performed component-wise on the characteristic variables.

\subsection{Extensions for Two-Dimensional Problems}
\label{subsec: twoD-scheme}

We briefly present, here, an extension of the WENO-RBF3 scheme to two-dimensional scalar conservation laws of the form
\begin{align}\label{eq2d:claw}
\begin{cases}
     & u_t + f(u)_x + g(u)_y = 0,\quad   \ (x,y)\in \RR \times \RR,\\
     &     u(x,y,0)    = u_0(x,y).
\end{cases}
\end{align}
Using the definition of the cell averages of the solution $u$
\begin{align}\label{eq2d:cell_avg}
\bar{u}_{ij} := \bar{u}_{ij}(t) :=\frac{1}{\Delta x \Delta y}\int_{x_{i-\frac{1}{2}}}^{x_{i+\frac{1}{2}}} \int_{y_{j-\frac{1}{2}}}^{y_{j+\frac{1}{2}}} u(x,y) \, dy \, dx,
\end{align}
we can recast equation \eqref{eq2d:claw} in its semi-discrete form, which is given by
\begin{equation}\label{eq2d:claw2}
	\frac{d}{dt} \bar{u}_{ij}(t) = - \frac{1}{\Delta x} \left( f_{i+\frac{1}{2},j} - f_{i-\frac{1}{2},j} \right) - \frac{1}{\Delta y} \left( g_{i,j+\frac{1}{2}} - g_{i,j-\frac{1}{2}} \right),
\end{equation}
where we have defined the fluxes
\begin{align}\label{eq2d:flux}
    \begin{split}
        f_{i\pm \frac{1}{2},j} &= \frac{1}{\Delta y}  \int_{y_{j-\frac{1}{2}}}^{y_{j+\frac{1}{2}}} f \left(u(x_{i\pm \frac{1}{2}}, y) \right) \, dy, \\
        g_{i,j\pm \frac{1}{2}} &= \frac{1}{\Delta x}  \int_{x_{i-\frac{1}{2}}}^{x_{i+\frac{1}{2}}} g \left(u(x,y_{j\pm \frac{1}{2}}) \right) \, dx.
    \end{split}
\end{align}
We discuss the construction of the flux $f_{i+ \frac{1}{2},j}$ in \eqref{eq2d:flux}, as the components $f_{i- \frac{1}{2},j}$ and $g_{i,j\pm \frac{1}{2}}$ can be obtained in an analogous manner. The integrals in the flux $f_{i+ \frac{1}{2},j}$ can be discretized with numerical quadrature, such as Gaussian quadratures or analytical integration using polynomial interpolation of a specified degree. For this work, we adopt the former approach, and employ Gauss-Legendre quadrature to perform the integration. In either case, with a selection of $N$ integration points $\{y_{\alpha} \in [y_{j-\frac{1}{2}},y_{j+\frac{1}{2}}]: \alpha=1, \cdots, N\}$ this leads to a discretization of the form
\begin{align}\label{eq2d:flux1}
    f_{i+ \frac{1}{2},j} \approx \frac{1}{\Delta y}\sum_{\alpha = 1}^N w_{\alpha} f \left(u(x_{i+ \frac{1}{2}},y_{\alpha}) \right),
\end{align}
with the corresponding integration weights $w_{\alpha}$ for $\alpha = 1, \cdots, N$.
The flux $f \left(u(x_{i+ \frac{1}{2}},y_{\alpha}) \right)$ appearing in \eqref{eq2d:flux1} will be replaced by the numerical flux
\begin{align}
    \hat{f} \left(u_{i+\frac{1}{2},\alpha}^{-}, u_{i+\frac{1}{2},\alpha}^{+} \right),
\end{align}
where we have used $u_{i+\frac{1}{2},\alpha}^{-}$ and $u_{i+\frac{1}{2},\alpha}^{+}$ to denote the left and right states, respectively, which are taken about the interface $x = x_{i+\frac{1}{2}}$, along the integration points $y_{\alpha}$ for $\alpha = 1, \cdots, N$.

\begin{figure}[htb!]
    \centering
    \includegraphics[width=0.7\textwidth, trim={5cm 10cm 5cm 5cm}, clip]{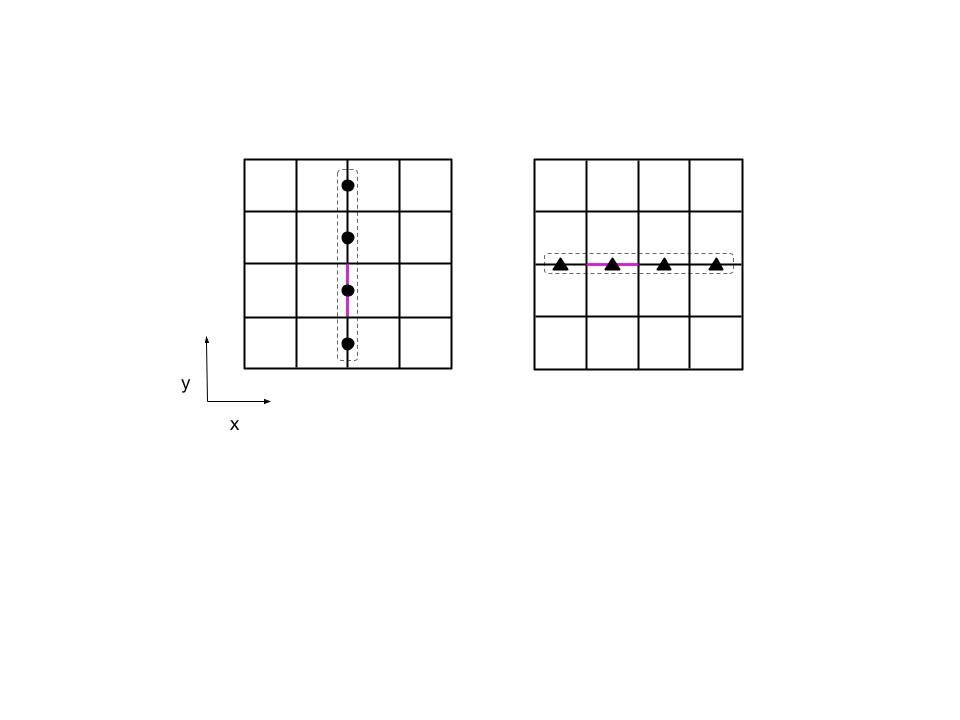}
    \caption{The $4 \times 4$ reconstruction patch used in the fourth order accurate WENO-RBF3 scheme. The first set of WENO reconstructions provides the data indicated by circles and triangles. This data is then used to construct interpolating functions that are projected onto the corresponding quadrature points, which lie on faces of cell $(i,j)$, i.e., $ u_{i+\frac{1}{2},\alpha}^{-}$ (left) and $ u_{\alpha, j+\frac{1}{2}}^{-}$ (right), shown in magenta. We wish to emphasize that the data on remaining faces of cell $(i,j)$ can be reconstructed using this \textit{same} $4 \times 4$ patch of data, since the approximations used in the shape parameter do not have to be properly winded.}
    \label{fig:2D_FV}
\end{figure}

Following Remark \ref{rmk:odd effective stencil}, the WENO-RBF3 scheme requires (for fourth order accuracy) a $4\times 4$ patch of cell averages $\bar{u}_{\ell,m}$, for $\ell = i-1, \cdots, i+2$ and $m = j-1, \cdots, j+2$,
to construct $u_{i+\frac{1}{2},\alpha}^{-}$. First, for each $m=j-1,\cdots,j+2$, using a set of 4 points $\{\bar{u}_{\ell, m}:\ell = i-1,\cdots,i+2 \}$, we construct the approximation $Av_m^{-}$ to the one-dimensional averages of the solution $v_m^{-}$ with respect to $y$ direction
\begin{align}
    v_m = \frac{1}{\Delta y} \int_{y_{m-\frac{1}{2}}}^{y_{m+\frac{1}{2}}} u \left(x_{i\pm \frac{1}{2}}, y \right) dy,
\end{align}
by applying the WENO-RBF3 scheme discussed in section \ref{sec:newscheme}.
Now that we have produced the 4 points $\{Av_m^{-}:m = j-1,\cdots,j+2 \}$, the point-wise values along the quadrature points can be obtained through interpolation. When the data is smooth, we can apply fixed stencil reconstructions at each of the integration points, e.g., 
\begin{align}\label{eq:Lagrange}
    u_{i+\frac{1}{2},\alpha}^{-} = \sum_{m=j-1}^{j+2} L(y_{\alpha}) Av_m^{-}, \quad \alpha = 1, \cdots, N,
\end{align}
where $L(y)$ is a Lagrange polynomial obtained from the points $\{Av_m^{-}\}$. In the event that the data is no longer guaranteed to be smooth, we can, instead, use a variant of the WENO-RBF3 scheme to perform the interpolation at quadrature points, rather than equation \eqref{eq:Lagrange}. An diagram of the patch used in the reconstructions for $ u_{i+\frac{1}{2},\alpha}^{-}$ and $ u_{\alpha, j+\frac{1}{2}}^{-}$ is shown in Figure \ref{fig:2D_FV}. If WENO schemes are used to perform interpolations at quadrature points, the interpolation coefficients on the substencils now depend on the evaluation point, in addition to certain information about the cell in which the reconstructions are performed. Consequently, this also holds for the linear WENO weights $\{ d_k \}$, which may cause (some of) the linear weights to become negative \cite{ShiHuShu}. For this reason, we plan to defer these quadrature reconstructions to future work. Instead, the developments provided in this work will employ fixed stencil reconstructions, such as \eqref{eq:Lagrange}, at the quadrature points.

\section{Numerical Results} \label{sec: Numerical results}

In this section, we provide some experimental results that demonstrate the improvements offered by the non-polynomial basis. We begin by investigating the performance of the proposed methods on some one-dimensional benchmark problems for the Euler equations before testing the method on the more challenging pressureless Euler system, which is weakly hyperbolic. All test problems use a third order explicit strong-stability preserving Runge-Kutta method \cite{gottlieb2001strong} for time integration. In our experiments, we compute the cell average values, from the initial conditions, using high-order Gauss-Legendre quadrature. We compare the proposed methods, which are labelled as WENO-RBF3 and WENO-RBF4, with third order and fifth order classical WENO schemes  \cite{Jiang1} and WENO-Z schemes \cite{bor}. These are labelled as WENO-JS3/WENO-Z3 and WENO-JS5/WENO-Z5, respectively. Unless otherwise stated, the results for the WENO-RBF4 scheme used a second order approximation (i.e., $p =2$) for the shape parameter \eqref{eq:lam_choice}. Lastly, note that the implementation of the WENO-RBF4 method uses the hybrid strategy discussed in section \ref{subsec: Hybrid scheme}. 

\subsection{Hyperbolic System}

We now present numerical results for the one-dimensional Euler equations of gas dynamics
\begin{equation}\label{euler1}
    U_t + F(U)_x = 0,
\end{equation}
where $V$ and $F(V)$ are given as
\begin{align*}
& U=(\rho , \rho u , E)^T,\hspace{.35cm} \\
& F(U)=(\rho u,\rho u^2+p,u(E+p))^T.
\end{align*}
Here $\rho, p, u$ and $E$ are the density, pressure,
velocity, and total energy, respectively. Additionally, we use the equation of state
\begin{equation*}
    p = (\gamma - 1)(E - \frac{1}{2}\rho u^2), \quad c = \sqrt{\frac{\gamma p}{\rho}},
\end{equation*}
with $\gamma = 1.4$. Here, $c$ denotes the local speed of sound in the gas.

All test problems for the hyperbolic system use the HLLC Riemann solver \cite{toro1} to compute the flux. Following \cite{hesthaven_HCL_book}, the HLLC flux is of the form
\begin{equation}
    h(U_l, U_r)^T = 
    \begin{cases}
        F(U_l), & \text{ if $s^{-} \geq 0$}, \\
        F_{l}^{*}, & \text{ if $s^{*} \geq 0 \geq s^{-}$}, \\
        F_{r}^{*}, & \text{ if $s^{+} \geq 0 \geq s^{*}$}, \\
        F(U_r), & \text{ if $s^{+} \leq 0$},
    \end{cases}
\end{equation}
where the intermediate velocity is defined as
\begin{equation*}
    s^{*} = \frac{ p_{r} - p_{l} + \rho_l u_l (s^{-} - u_l) - \rho_r u_r (s^{+} - u_r)  }{ \rho_l (s^{-} - u_l ) - \rho_r (s^{+} - u_r ) }.
\end{equation*}
Estimates for the minimum and maximum wave speeds given, respectively, by
\begin{equation*}
    s^{-} = \min \left( u_l - c_l, u_l, u_l + c_l \right), \quad s^{+} = \min \left( u_r - c_r, u_r, u_r + c_r \right).
\end{equation*}
Using the Roe pressure
\begin{equation*}
    p_{lr} = \frac{p_r + p_l + \rho_l(s^{-} - u_l)(s^{*} - u_l) + \rho_r(s^{*} - u_r)(s^{+} - u_r) }{2},
\end{equation*}
we form the intermediate fluxes $F_{l}^{*}$ and $F_{r}^{*}$
\begin{equation*}
    F_{l}^{*} = \frac{s^{*}\Big( s^{-}u_l - F(U_l) \Big) + s^{-} p_{lr} \left(0, 1, s^{*}\right)^{T} }{s^{-} - s^{*}}, \quad F_{r}^{*} = \frac{s^{*}\Big( s^{+}u_r - F(U_r) \Big) + s^{-} p_{lr} \left(0, 1, s^{*}\right)^{T} }{s^{+} - s^{*}}.
\end{equation*}


\begin{example}\label{ex:smooth}
{\rm
First, we consider the smooth advection problem. Using the initial data
\begin{equation*}
    \rho(x) = 1 + 0.5\sin(4\pi x),
\end{equation*}
along with the selections $u = 1$, and $p = 1$, the system reduces to a single advection equation, subject to periodic boundary conditions. The exact solution on $[0,1]$ is given by
\begin{equation*}
    \rho(x,t) = 1 + 0.5\sin \Big(4\pi \left(x - ut \right) \Big).
\end{equation*}
In Table \ref{tab:euler},  we present both the numerical errors and convergence orders of the proposed scheme. When setting the shape parameter $\lambda^2$ for the RBF $\phi=e^{-\lambda^2x^2}$, the approximations in \eqref{eq:lam_choice} were computed to first and second order accuracy (i.e., $p=1$ and $2$ in \eqref{eq:lam_choice}). These choices are reflected under the respective labels ``Proposed I" and ``Proposed II" shown in Table \ref{tab:euler}. The results demonstrate the expected improved convergence through the choice of the shape parameter. With regard to the three-point methods, the proposed RBF scheme achieves its intended accuracy along with a noticeable reduction in the errors compared to alternative reconstruction methods.  
}

\begin{table}[t!]
	\centering
	
	\begin{tabular}{r|c|ccc}
		\hline\hline
	\multirow{5}{*} {$L_\infty$} &	$N_x$ &   WENO-JS3  &  WENO-Z3 &  RBF3 (Proposed)    \\
		 \cline{2-5}
        & 20 & 3.04e-01 ( \,\,---\,\, ) & 2.78e-01 ( \,\,---\,\, ) & 1.66e-02 ( \,\,---\,\, ) \\
        & 40 & 1.57e-01 ( 0.96 ) & 8.59e-02 ( 1.69 ) & 6.01e-04 ( 4.79 ) \\
        & 80 & 4.81e-02 ( 1.71 ) & 2.78e-02 ( 1.63 ) & 2.06e-05 ( 4.87 ) \\
        &160 & 1.58e-02 ( 1.61 ) & 4.66e-03 ( 2.58 ) & 8.06e-07 ( 4.68 ) \\
        &320 & 2.80e-03 ( 2.50 ) & 4.45e-04 ( 3.39 ) & 7.75e-08 ( 3.38 ) \\
		\hline\hline
	\multirow{5}{*} {$L_1$} &	$N_x$  &   WENO-JS3  &  WENO-Z3 &  RBF3 (Proposed)    \\
		 \cline{2-5}
        &  20 & 4.73e-01 ( \,\,---\,\, ) & 4.38e-01 ( \,\,---\,\, ) & 2.67e-02 ( \,\,---\,\, ) \\
        &  40 & 2.79e-01 ( 0.76 ) & 1.87e-01 ( 1.23 ) & 1.08e-03 ( 4.63 ) \\
        &  80 & 1.11e-01 ( 1.33 ) & 6.38e-02 ( 1.55 ) & 3.60e-05 ( 4.90 ) \\
        & 160 & 4.19e-02 ( 1.41 ) & 1.78e-02 ( 1.85 ) & 1.28e-06 ( 4.82 ) \\
        & 320 & 1.22e-02 ( 1.78 ) & 3.09e-03 ( 2.52 ) & 1.25e-07 ( 3.35 ) \\
		\hline\hline

		\hline\hline
	\multirow{5}{*} {$L_\infty$} &	$N_x$ &   WENO-JS5  &  RBF4 (Proposed I) &     RBF4 (Proposed II)   \\
		 \cline{2-5}
		& 20  & 1.17e-02 ( \,\,---\,\, )&	 3.12e-03 ( \,\,---\,\, )&	 3.03e-04 ( \,\,---\,\, ) \\
		& 40  & 6.95e-04 ( 4.07 ) &	 1.02e-04 ( 4.94 ) &	 5.70e-06 ( 5.73 )	 \\
		& 80  & 2.61e-05 ( 4.74 ) &	 3.21e-06 ( 4.99 ) &	 8.71e-08 ( 6.03 )	 \\
		&  160& 8.70e-07 ( 4.91 ) &	 1.00e-07 ( 5.00 ) &	 1.39e-09 ( 5.97 )	 \\
		&  320& 2.68e-08 ( 5.02 ) &	 3.14e-09 ( 5.00 ) &	 2.16e-11 ( 6.01 )	 \\
		\hline\hline
	\multirow{5}{*} {$L_1$} &	$N_x$  &   WENO-JS5  &  RBF4 (Proposed I) &     RBF4 (Proposed II)   \\
		 \cline{2-5}
        & 20 & 8.14e-03 ( \,\,---\,\, )&	 1.60e-03 ( \,\,---\,\, )&	 1.61e-04 ( \,\,---\,\, ) \\
        &   40 & 3.74e-04 ( 4.45 ) &	 5.29e-05 ( 4.92 ) &	 3.31e-06 ( 5.61 ) \\
        &   80 & 1.18e-05 ( 4.98 ) &	 1.68e-06 ( 4.98 ) &	 5.31e-08 ( 5.96 ) \\
        &  160 & 3.69e-07 ( 5.00 ) &	 5.27e-08 ( 5.00 ) &	 8.39e-10 ( 5.98 ) \\
        &  320 & 1.15e-08 ( 5.00 ) &	 1.65e-09 ( 5.00 ) &	 1.31e-11 ( 6.00 ) \\
		\hline\hline
	\end{tabular}
\medskip
\caption{(Example \ref{ex:smooth}) $L_\infty $ and $L_1$ errors and convergence rates (*) at $t=1$.}
	\label{tab:euler}
\end{table}
\end{example}

\begin{example}\label{ex:lax}
{\rm  We test our method  using the Lax problem \cite{lax1} with the initial condition
\begin{equation*}
	(\rho, u, p) =
	\begin{cases}
	(0.445,0.698,3.528) & \text{if $x\in [ -5,  0)$}, \\
	(0.50, 0, 0.571) & \text{if   $x\in [0,  5]$}.
	\end{cases}
\end{equation*}
The numerical results for the density profiles are displayed in Figure \ref{fig:lax}, at time $t = 1.3$, using $\Delta x = 10/200$. 

Among three-point WENO schemes, the proposed method demonstrates improvement in capturing the jumps in the profile, when compared to other three-point WENO methods. In the case of the four-point RBF schemes, the method provides improvements over WENO-JS and is comparable to the solution generated with WENO-Z5. This similarity may be a consequence of the size of the substencils employed by each of the methods. The substencils used in WENO-JS5 and WENO-Z5 consist of three points, whereas the proposed four-point RBF methods use substencils with two points. This discrepancy in the behavior of the method is something we wish to address in a subsequent paper.
}

\begin{figure}[htb!]
    \centering
        \subfigure[Comparison of three-point WENO schemes.]{  
    \includegraphics[width=0.6\textwidth]{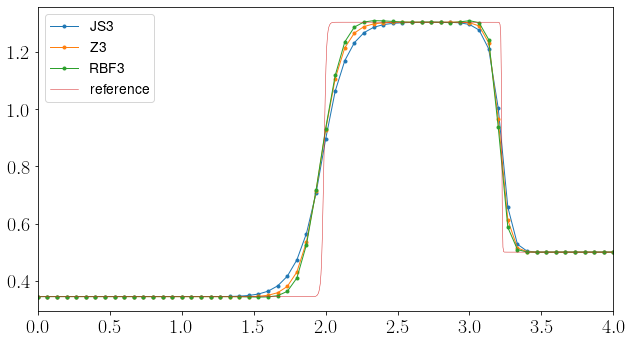}
    }
        \subfigure[Comparison with five-point WENO schemes.]{  
    \includegraphics[width=0.6\textwidth]{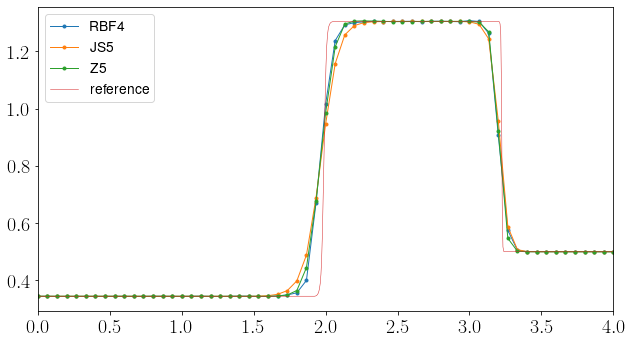}
    }
    \caption{(Example \ref{ex:lax}) Density profiles of Lax problem \cite{lax1} at $t=1.3$ with $\Delta x = 10/150$.}
    \label{fig:lax}
\end{figure}

\end{example}

\begin{example}\label{ex:sod}
{\rm As a next example, let us consider the Sod problem \cite{sod}
with the initial condition
\begin{equation*}
	(\rho, u, p) =
	\begin{cases}
	(1.000,0.750,1.000) & \text{if $x\in [0,  0.5)$}, \\
	(0.125,0.000,0.100) & \text{if $x\in [0.5, 1]$}.
	\end{cases}
\end{equation*}
The numerical results for density profiles are given in Figure \ref{fig:sod}
at time $t = 0.2$ using $\Delta x = 1/100$.  Here, we observe results which are similar to the Lax problem (see Example \ref{ex:lax}). As noted above, we plan to investigate these shortcomings in our future work involving WENO-RBF methods.
}

\begin{figure}[htb!]
    \centering
        \subfigure[Comparison of three-point WENO schemes.]{  
    \includegraphics[width=0.6\textwidth]{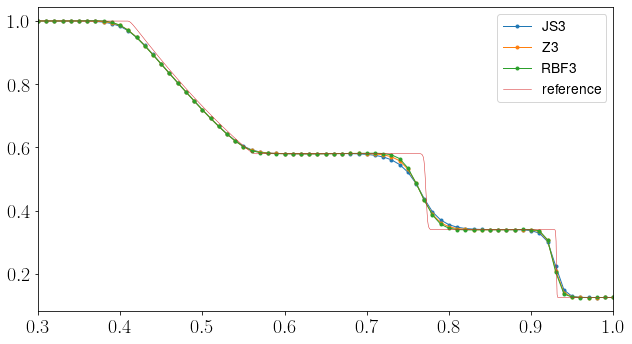}
    }
        \subfigure[Comparison with five-point WENO schemes.]{  
    \includegraphics[width=0.6\textwidth]{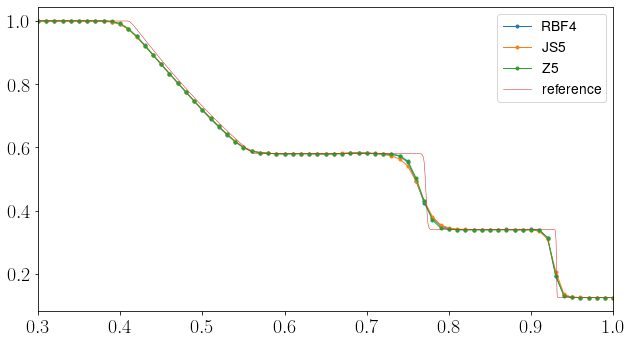}
    }
    \caption{(Example \ref{ex:sod}) Density profiles of Sod problem \cite{sod} at $t=0.2$ using $\Delta x = 1/100$.}
    \label{fig:sod}
\end{figure}

\end{example}

\begin{example}\label{ex:shu}
{\rm We now look at the Shu-Osher problem \cite{shu2}, which uses Riemann initial data for the shock entropy wave
interaction.
The approximate solutions are computed on the interval $[-5,5]$
with the initial state
\begin{equation}\label{eq:shu}
(\rho, u, p) =
\begin{cases}
(3.857143,2.629369,10.33333) & \text{if $x\in [-5, -4)$}, \\
(1+ \varepsilon \sin(kx),0,1) & \text{if $x\in [-4,   5]$},
\end{cases}
\end{equation}
where $k$ and $\varepsilon$ denote the wave number and amplitude
of the entropy wave respectively.
We take $k=5$ and $\varepsilon = 0.2$ in our experiments.

In Figure \ref{fig:shu}, we provide plots which compare the proposed WENO-RBF methods against three-point and five-point WENO methods. We observe the largest improvement in the comparison of three-point WENO schemes, with the WENO-RBF3 method producing results with less dissipation than WENO-JS3 and WENO-Z3. With regard to the four-point method, i.e., WENO-RBF4, improvements over WENO-JS5 are quite clear. These improvements are less apparent when compared against WENO-Z5, which seems to perform better in certain areas. As discussed in Example \ref{ex:lax}, part of this improvement may be attributed to the larger substencil size used by the WENO-Z5 method.

}

\begin{figure}[htb!]

    \begin{center}
        \subfigure[Comparison of three-point WENO schemes.]{  
        \includegraphics[width=0.6\textwidth]{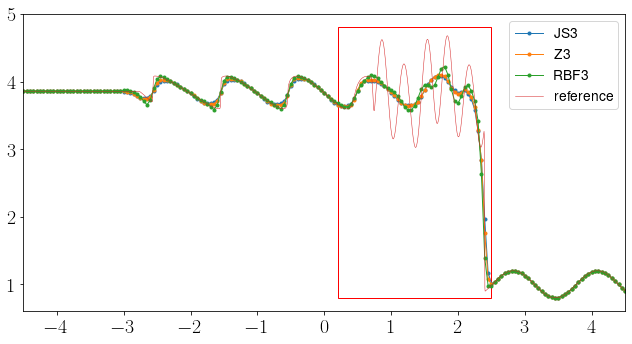}
        \includegraphics[width=0.33\textwidth]{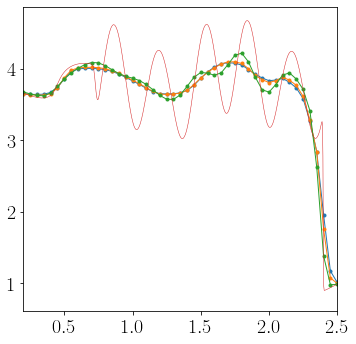}
        }
        
        \subfigure[Comparison with five-point WENO schemes.]{  
        \includegraphics[width=0.6\textwidth]{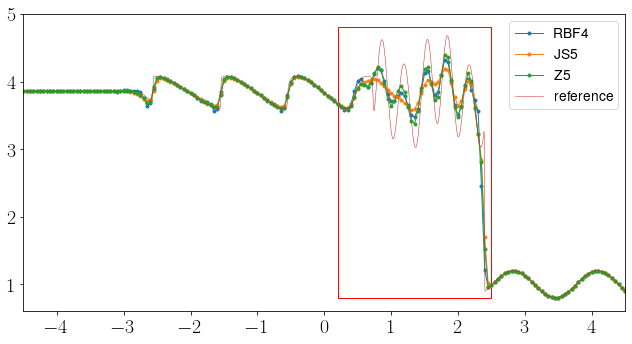}
        \includegraphics[width=0.33\textwidth]{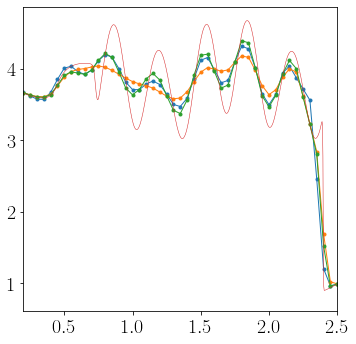}
        }
        
    \end{center}
    \caption{(Example \ref{ex:shu}) Density profiles of the shock-entropy interaction problem at $t=1.8$, with $\Delta x = 10/300$.}
    \label{fig:shu}
\end{figure}

\end{example}

\begin{example}\label{ex:toro}

{\rm 
In \cite{toro}, Titarev and Toro suggested the following Riemann initial data for the shock entropy wave interaction:
\begin{equation*}
(\rho, u, p) =
\begin{cases}
( 1.515695,0.523346,1.80500) & \text{if $x\in [-5, -4.5)$}, \\
(1+ 0.1 \sin(20\pi x),0,1) & \text{if $x\in [-4.5,   5]$}.
\end{cases}
\end{equation*}
on the interval $[-5,5]$.
Figure \ref{fig:toro} shows the results for this test case at time at $t=5$. This problem allows us to test the method in environments that support highly oscillatory structures. Compared to the Shu-Osher problem (see Example \ref{ex:shu}), the prescribed initial density exhibits a much larger wave number, i.e., $20$. 

Among the three-point WENO schemes, we see from Figure \ref{fig:toro_a} that the WENO-RBF3 method offers a clear improvement over WENO-JS3 and WENO-Z3. We also find that WENO-RB4 yields an improvement over the results obtained for the Shu-Osher problem, with regard to capturing complex wave patterns. While there is a clear improvement in capturing the features of the small wavelength oscillations, we do note the slight undershoot and overshoot in the vicinity of $x=-1.9$, where the transition into the oscillatory region occurs. This behavior is far less apparent in the WENO-JS5 and WENO-Z5 methods.

}

\begin{figure}[htb!]
    \begin{center}
     
        \includegraphics[width=0.86\textwidth]{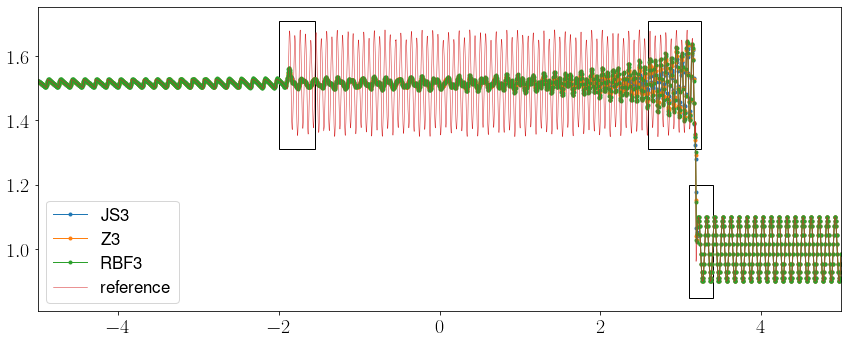}
            \subfigure[Comparison of three-point WENO schemes.]{ 
        \includegraphics[width=0.29\textwidth]{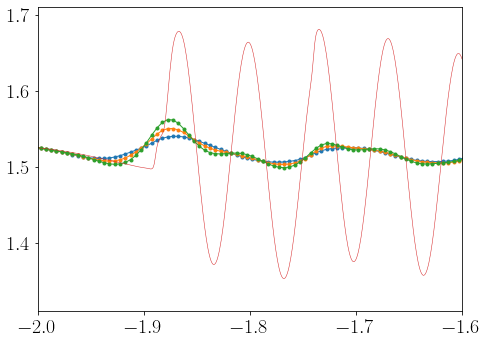}
        \includegraphics[width=0.336\textwidth]{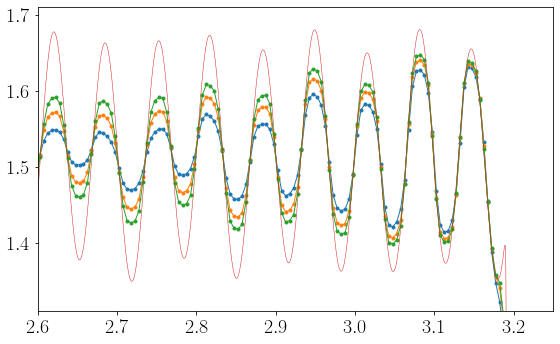}
        \includegraphics[width=0.218\textwidth]{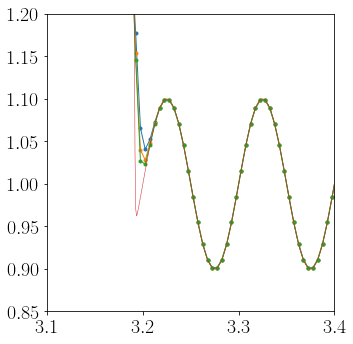}
        \label{fig:toro_a}
        }
        
        \includegraphics[width=0.86\textwidth]{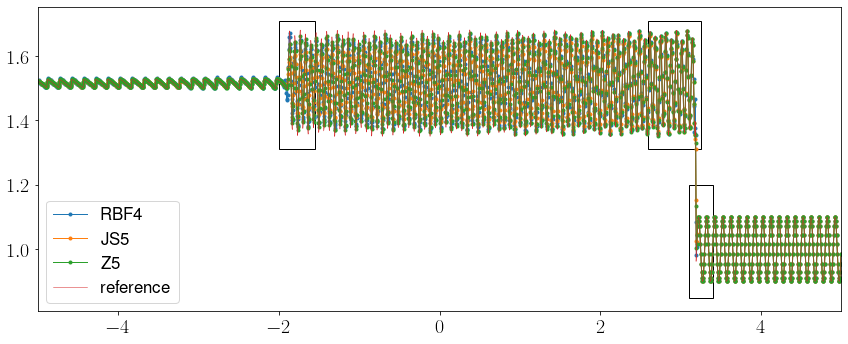}
            \subfigure[Comparison with five-point WENO schemes.]{ 
        \includegraphics[width=0.29\textwidth]{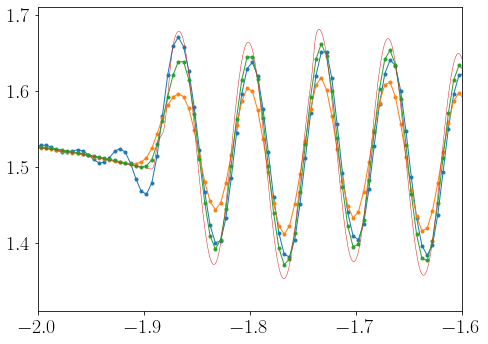}
        \includegraphics[width=0.336\textwidth]{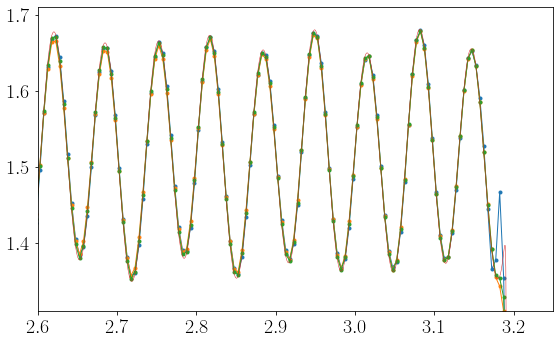}
        \includegraphics[width=0.218\textwidth]{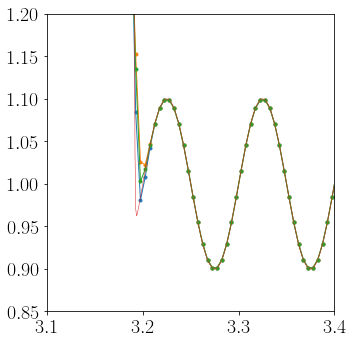}
        \label{fig:toro_b}
        }
    \end{center}
    
    \caption{(Example \ref{ex:toro}) Density profiles of the shock-turbulence problem \cite{toro} at $t=5$ using $\Delta x=1/200$.}
    \label{fig:toro}
\end{figure}

\end{example}


\subsection{\bf Weakly Hyperbolic System}
\label{subsec: Weakly hyperbolic system}

In this subsection, we consider the more challenging pressureless Euler system, which, in one space dimension, is given by
\begin{equation}\label{euler_p0}
U_t + F(U)_x = 0,
\end{equation}
with
\begin{equation*}
U=(\rho , \rho u )^T,\hspace{.35cm} F(U)=(\rho u,\rho
u^2)^T.
\end{equation*}
The system \eqref{euler_p0} is important in modeling systems of dilute gases in a vacuum which undergo few collision events. When collisions do occur, they are said to be perfectly inelastic, which causes the gas particles to stick together. These collisions lead to the emergence of so-called $\delta$-shocks, which are the primary feature of interest in these models. Note that because the system is weakly hyperbolic, the characteristic decomposition is not available for this problem. Consequently, the simple Lax-Friedrichs flux cannot be directly used in this problem. Instead, we consider the Godunov flux, outlined in \cite{YWS}, which was originally introduced by Bouchut, \etal \cite{Bouchut_pressureless_euler}. Before defining the flux, we remark that the latter work contains a fairly diverse collection of literature on isothermal gas dynamics, so we refer the interested reader to references therein for further details. The former article considered more general problems involving $\delta$-shocks, including the system \eqref{euler_p0}, and applied DG methods to solve such problems. To define the flux, suppose we have left and right numerical approximations $U_{l}=(\rho _l,\rho_l u_l)^T$ and $U_{r}=(\rho _r,\rho_r u_r)^T$. Then, the numerical flux is given by
\begin{equation*} 
h(U_l, U_r)^T = \begin{cases}
    \left(\rho_l u_l, \rho_l u^2_l \right)^{T},  & \text{ if $u_l>0, u_r>0$}, \\
    (0,0)^{T}, & \text{ if $u_l\leq 0, u_r>0$}, \\
    \left(\rho_r u_r, \rho_r u^2_r \right)^{T}, & \text{ if $u_l\leq 0, u_r\leq 0$}, \\
    \left( \rho_l u_l, \rho_l u^2_l \right)^{T}, & \text{ if $u_l>0, u_r\leq 0, v>0$}, \\
    \left( \rho_r u_r, \rho_r u^2_r \right)^{T}, & \text{ if $u_l>0, u_r\leq 0, v<0$}, \\
    \left( \dfrac{\rho_l u_l+\rho_r u_r}{2}, \rho_l u^2_l=\rho_r u^2_r \right)^{T} , & \text{ if $u_l>0, u_r\leq 0, v=0$},
\end{cases}
\end{equation*}
where $$ v = \frac{\sqrt{\rho_l}u_l+\sqrt{\rho_r}u_r}{\sqrt{\rho_l}+\sqrt{\rho_r}}. $$

\begin{example} 
\label{euler_p_smooth}
{\rm 
        For our first test, we seek to determine the approximation order for the weak hyperbolic system. To this end,
		we solve pressureless Euler system \eqref{euler_p0} with the following initial data
		$$ \rho _0 (x)=\sin(x)+2, \quad u_0(x)=\sin(x)+2,$$
		subject to periodic boundary conditions. The exact solution for this problem is
		$$ \rho (x,t)=\frac{\rho _0(x_0)}{1+t u'_0(x_0)}, \quad u(x,t)=u_0(x_0),  $$
		where $x_0$ is given implicitly by
		$$ x_0 + tu_0(x_0)=x.$$
		The $L_{\infty}$ and $L_1$ errors and approximation orders for the density $\rho$ at $t=0.1$ are given in Table \ref{tab:peuler}. For the three-point methods, we find that both WENO-JS3 and WENO-Z3 exhibit difficulties in achieving second order accuracy. While some reduction in the error is attained by switching from WENO-JS3 to WENO-Z3, we find that convergence properties are fairly similar between these two methods. In contrast, the WENO-RBF3 method achieves the intended convergence rate along with significant reduction in the error. 
		
		Similarly, in the case of the four-point method (WENO-RBF4), we achieve the intended accuracy of the method. In accordance with Theorem \ref{thm:even}, the shape parameters for the RBF $\phi=e^{-\lambda^2x^2}$, were computed to first and second order accuracy (i.e., $p=1$ and $p=2$ in \eqref{eq:lam_choice}) and have been assigned the corresponding labels ``Proposed I" and ``Proposed II", as in the Example \ref{ex:smooth}. Moreover, these schemes exhibit convergence rates, of fifth and sixth order accuracy, respectively, as shown in Table \ref{tab:peuler}. In contrast, the WENO-JS5 scheme suffers a reduction in the convergence order by nearly a factor of two, a feature which was also observed in \cite{YWS}. 
}
\end{example}

\begin{table}[htb!]
	\centering
	
	\begin{tabular}{r|c|ccc}
		\hline\hline
	\multirow{5}{*} {$L_\infty$} &	$N_x$  &   WENO-JS3  &  WENO-Z3 &  RBF3 (Proposed)    \\
		 \cline{2-5}
        &  20 & 6.27e-02 ( \,\,---\,\, ) & 6.16e-02 ( \,\,---\,\, ) & 7.87e-03 ( \,\,---\,\, ) \\
        &  40 & 6.20e-02 ( 0.02 ) & 4.86e-02 ( 0.34 ) & 3.96e-03 ( 0.99 ) \\
        &  80 & 1.26e-01 (-1.03 ) & 5.82e-02 (-0.26 ) & 3.19e-04 ( 3.64 ) \\
        & 160 & 3.15e-02 ( 2.00 ) & 1.68e-02 ( 1.79 ) & 2.06e-05 ( 3.95 ) \\
        & 320 & 1.66e-02 ( 0.93 ) & 5.50e-03 ( 1.62 ) & 1.22e-06 ( 4.08 ) \\
		\hline\hline
	\multirow{5}{*} {$L_1$} &	$N_x$    &   WENO-JS3  &  WENO-Z3 &  RBF3 (Proposed)    \\
		 \cline{2-5}
        &  20 & 8.88e-02 ( \,\,---\,\, ) & 8.16e-02 ( \,\,---\,\, ) & 1.64e-02 ( \,\,---\,\, ) \\
        &  40 & 3.97e-02 ( 1.16 ) & 2.83e-02 ( 1.53 ) & 6.91e-03 ( 1.25 ) \\
        &  80 & 4.46e-02 (-0.16 ) & 1.42e-02 ( 0.99 ) & 5.37e-04 ( 3.69 ) \\
        & 160 & 9.68e-03 ( 2.20 ) & 5.41e-03 ( 1.39 ) & 3.39e-05 ( 3.98 ) \\
        & 320 & 2.81e-03 ( 1.78 ) & 8.37e-04 ( 2.69 ) & 2.21e-06 ( 3.94 ) \\
		\hline\hline

		\hline\hline
	\multirow{5}{*} {$L_\infty$} &	$N_x$  &   WENO-JS5  &  RBF4 (Proposed I) &     RBF4 (Proposed II)    \\
		 \cline{2-5}
		& 20  & 1.89e-03 ( \,\,---\,\, )&	 8.38e-04 ( \,\,---\,\, ) &	 7.38e-04 ( \,\,---\,\, ) \\
		& 40  & 6.90e-05 ( 4.77 ) &	 1.98e-05 ( 5.41 ) &	 1.32e-05 ( 5.80 )	 \\
		& 80  & 3.62e-06 ( 4.25 ) &	 5.68e-07 ( 5.12 ) &	 2.12e-07 ( 5.96 )	 \\
		&  160& 2.48e-07 ( 3.87 ) &	 1.67e-08 ( 5.08 ) &	 3.37e-09 ( 5.97 )	 \\
		&  320& 1.90e-08 ( 3.71 ) &	 5.05e-10 ( 5.05 ) &	 5.30e-11 ( 5.99 )	 \\
		\hline\hline
	\multirow{5}{*} {$L_1$} &	$N_x$    &   WENO-JS5  &  RBF4 (Proposed I)  &     RBF4 (Proposed II)    \\
		 \cline{2-5}
        & 20 & 3.23e-03 ( \,\,---\,\, ) &	 1.95e-03 ( \,\,---\,\, )&	 1.45e-03 ( \,\,---\,\, ) \\
        & 40 & 1.13e-04 ( 4.84 ) &	 4.52e-05 ( 5.43 ) &	 2.41e-05 ( 5.91 ) \\
        & 80 & 3.68e-06 ( 4.94 ) &	 1.15e-06 ( 5.30 ) &	 3.81e-07 ( 5.98 ) \\
        &160 & 1.25e-07 ( 4.88 ) &	 3.17e-08 ( 5.17 ) &	 6.04e-09 ( 5.98 ) \\
        &320 & 4.68e-09 ( 4.74 ) &	 9.27e-10 ( 5.10 ) &	 9.47e-11 ( 5.99 ) \\
		\hline\hline
	\end{tabular}

\medskip
\caption{(Example \ref{euler_p_smooth}) $L_\infty $ and $L_1$ errors and convergence rates (*) at $t=0.1$.}
	\label{tab:peuler}
\end{table}

\begin{example}\label{ex:dshock}
{\rm
We solve the problem known as two interactive blast wave problem \cite{wood}
which has the initial data
\begin{equation*}
	(\rho_0,u_0) =  \begin{cases}
		(1,1) & \text{if $ x <0 $}, \\
		(0.25,0) & \text{if $ x>0 $}.
	\end{cases}
\end{equation*}
In Figure \ref{fig:dshock}, we plot the density profiles for each of the methods at time $t = 0.3$.

Our results indicate that the three-point WENO methods provide fairly similar estimates of the shock width, but produce remarkable differences in the overall height of the $\delta-$shock, as shown in Figure \ref{fig:blast_a}. Of these methods, the shock generated by WENO-RBF3 method exhibits the sharpest peak, compared to WENO-JS3 and WENO-Z3. Similar observations can be made regarding the four-point WENO-RBF4 method, which is presented in Figure \ref{fig:blast_b}. The shock width predicted by the WENO-RBF4 method is, again, sharper and slightly more narrow than the that predicted by WENO-JS5 and WENO-Z5. To the right of the shock, we also observe some undershoots in the densities predicted by each of these methods; however, the undershoot in the RBF method is marginally smaller than WENO-JS5 and noticeably smaller than WENO-Z5. As discussed earlier, reducing the overshoot and undershoot in these rapidly varying transition regions is something we plan to address in our future work.

}

\begin{figure}[htb!]

    \begin{center}
        \subfigure[Comparison of three-point WENO schemes using $\Delta x = 1/200$.]{  
        \includegraphics[width=0.6\textwidth]{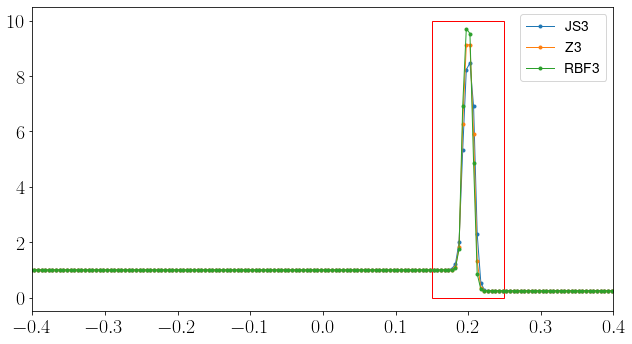}
        \includegraphics[width=0.33\textwidth]{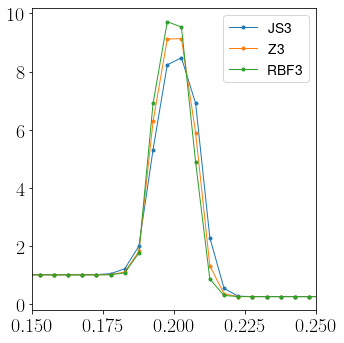}
        \label{fig:blast_a}
    }
    \subfigure[Comparison with five-point WENO schemes using $\Delta x = 1/80$.]{
        \includegraphics[width=0.6\textwidth]{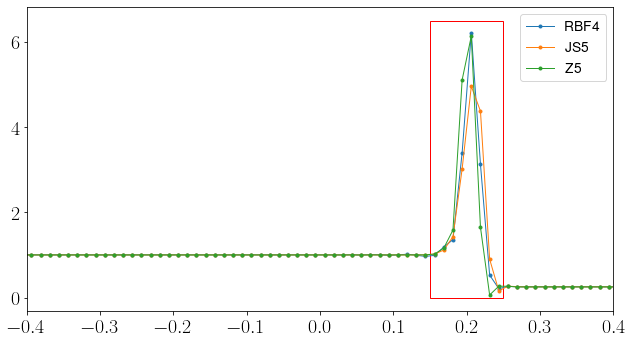}
        \includegraphics[width=0.33\textwidth]{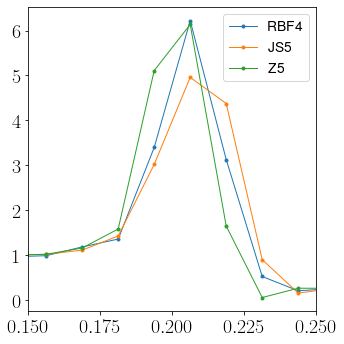}
        \label{fig:blast_b}
    }
    \end{center}
    \caption{(Example \ref{ex:dshock}) Density profiles of $\delta$-shock wave problem \cite{wood} at $t=0.3$ }
    \label{fig:dshock}
\end{figure}

\end{example}

\subsection{Two-Dimensional Scalar Hyperbolic Problem}

In this section, we present convergence results for a two-dimensional scalar problem using the extensions described in section \ref{subsec: twoD-scheme}. As discussed earlier, an extension to multi-dimensional problems can be achieved with line-by-line applications of the proposed one-dimensional WENO schemes, introduced earlier in this work (see e.g., section \ref{sec: WENO schemes}). We only provide results for a single two-dimension test problem, as a proof of concept, since this is not the central theme of this paper. In a subsequent article, we plan to focus our efforts on two and three-dimensional applications.

\begin{example}\label{ex:2dburgers}

{\rm 
In this example, we apply the three-point RBF scheme to a multi-dimensional problem, focusing, in particular, on the two-dimensional Burgers' equation
\begin{equation*}
    \label{eq:2D Burgers}
    u_t + \frac{1}{2} \left(u^2 \right)_x + \frac{1}{2} \left(u^2 \right)_y = 0.
\end{equation*}
We use the initial condition 
$$ u(x,0) = 0.5 + \sin(x+y), $$
along with periodic boundary conditions on the domain $[-\pi,\pi] \times [-\pi,\pi]$. For this problem, we used the local Lax-Friedrichs flux
\begin{equation}
    h(u^{-}, u^{+}) = \frac{f(u^{-}) + f(u^{+})}{2} - \alpha \frac{ u^{+} - u^{-} }{2},
\end{equation}
where $\alpha$ is the maximum wave speed taken over the local states
\begin{equation*}
    \alpha = \max \left( \left\lvert u^{-} \right\rvert, \left\lvert u^{+} \right\rvert \right).
\end{equation*}
Convergence of the smooth solutions was measured using the final time $T=0.2$ and we chose the timestep size according to $\Delta t = \min (\Delta x^{4/3}, \Delta y^{4/3} )$. We report the $L_{\infty}$ and $L_1$ errors and convergence rates in Table \ref{tab:2dburgers} for each of the three-point WENO schemes. In each of these methods, we used the same fifth order Lagrange interpolating polynomial \eqref{eq:Lagrange} to perform the reconstructions at the quadrature points.

Our preliminary results indicate that the WENO-RB3 method achieves, at least, third order accuracy and with the convergence rate tending to fourth order as the mesh resolution increases. Moreover, the errors in the proposed scheme, when compared to WENO-JS3 and WENO-Z3, are noticeably smaller. Given that each of these methods used the same quadrature reconstruction procedure, this improvement stems from the use of the new reconstruction method. We note that the convergence rates for the WENO-RBF3 method, in this example, are not ideal, which is something we are currently working to resolve.

}

\begin{table}[t]
	\centering
	
	\begin{tabular}{r|c|ccc}
		\hline\hline
	\multirow{5}{*} {$L_\infty$} &	$N_x \times N_y$ &   WENO-JS3 &  WENO-Z3  &  RBF3    \\
		 \cline{2-5}
        & 10 $\times$  10 &  8.2958e-02 ( \,\,---\,\, ) & 6.9456e-02 ( \,\,---\,\, )    & 1.9902e-02 ( \,\,---\,\, ) \\
        & 20 $\times$  20 &  4.5343e-02 ( 0.87 )        & 3.7121e-02 ( 0.90 )           & 1.1847e-02 ( 0.75 ) \\
        & 40 $\times$  40 &  1.9416e-02 ( 1.22 )        & 1.4812e-02 ( 1.33 )           & 2.2575e-03 ( 2.39 ) \\
        & 80 $\times$  80 &  7.7397e-03 ( 1.33 )        & 4.9918e-03 ( 1.57 )           & 2.6275e-04 ( 3.10 ) \\
        &160 $\times$ 160 &  2.5147e-03 ( 1.62 )        & 9.7351e-04 ( 2.36 )           & 2.0994e-05 ( 3.65 ) \\
		\hline\hline
	\multirow{5}{*} {$L_1$} &	$N_x \times N_y$  &   WENO-JS3   &  WENO-Z3  &  RBF3    \\
		 \cline{2-5}
        &  10 $\times$  10 & 1.2567e+00 ( \,\,---\,\, ) & 1.0833e+00 ( \,\,---\,\, )& 3.0803e-01 ( \,\,---\,\, ) \\
        &  20 $\times$  20 & 4.2019e-01 ( 1.58 )        & 2.8167e-01 ( 1.94 )       & 9.6608e-02 ( 1.67 ) \\
        &  40 $\times$  40 & 1.2268e-01 ( 1.78 )        & 7.9952e-02 ( 1.82 )       & 1.0978e-02 ( 3.14 ) \\
        &  80 $\times$  80 & 3.0699e-02 ( 2.00 )        & 1.6416e-02 ( 2.28 )       & 8.8519e-04 ( 3.63 ) \\
        & 160 $\times$ 160 & 6.3137e-03 ( 2.28 )        & 1.9687e-03 ( 3.06 )       & 6.9829e-05 ( 3.66 ) \\
		\hline\hline
	\end{tabular}
\medskip
\caption{(Example \ref{ex:2dburgers}) $L_\infty $ and $L_1$ errors and convergence rates (*) at $t=0.2$.}
\label{tab:2dburgers}
\end{table}

\end{example}

	\section{Conclusion} \label{sec: Conclusion}

In this paper, we proposed several FV WENO-RBF methods which achieve superconvergence using a non-polynomial basis consisting of RBFs. Superconvergence in the interpolation component of the reconstruction was achieved by exploiting the shape parameter in the definition of the RBF. More specifically, by analyzing the error in the reconstructions, we derived expressions for optimal shape parameters which improved the convergence order of the interpolation. While the methods developed in this work considered Gaussian functions, the same techniques can be easily applied (or adapted) to develop schemes based on other RBFs, as well as other bases. The proposed schemes, whose constructions make use of fairly compact stencils, incorporate new smoothness indicators which were previously shown to be highly effective at discerning rapid changes in a function, even on a small data stencil. To alleviate the heavy computational cost typically associated with WENO methods, we also implemented a hybrid solver that dynamically prescribes the reconstruction method according to the local smoothness of the function. The proposed schemes were compared with several competing methods using one-dimensional systems of conservation laws, along with a two-dimensional test problem to demonstrate extensions for multiple dimensions. While the proposed methods demonstrated improved shock-capturing capabilities, the use of the non-polynomial basis was shown to be particularly effective for problems exhibiting rapid transitions as well as complex wave structures and singularities. Furthermore, in the case of the pressureless Euler equations, the proposed WENO-RBF methods were experimentally shown to achieve their theoretical convergence rates, avoiding the order reduction encountered by competitive reconstruction techniques. A highlight of this work is reflected in the blast wave problem (see Example \ref{ex:dshock} in section \ref{subsec: Weakly hyperbolic system}), where the proposed RBF methods lead to markedly different predictions of the shock. While these results suggest several avenues for future research, we plan to conduct additional experiments to develop multi-dimensional algorithms and investigate strategies for further reducing oscillations surrounding transition regions and singularities.

\section*{Acknowledgements}

We would like to thank AFOSR and NSF for their support of this work under grants FA9550-19-1-0281 and FA9550-17-1-0394 and NSF grant DMS 191218.

	\begin{appendices}

\section{Approximating the shape parameter $\lambda^2$} \label{app:lambda}

We briefly present, here, how to compute the shape parameter $\lambda^2$ in the radial basis function $\phi (x) = \exp (-{\lambda}^2 x^2)$. As an example, we consider the (fixed stencil) RBF3 scheme with the shape parameter defined in \eqref{eq:lam4rbf3}, and we assume a set of cell average values $\{\bar{u}_k\}$ is available.

To construct the values $u_{j+1/2}^{-}$ and $u_{j+1/2}^{+}$, the reconstruction stencils are given, respectively, as $\{\bar{u}_{j-1},\bar{u}_j,\bar{u}_{j+1}\}$ and $\{\bar{u}_j,\bar{u}_{j+1}, \bar{u}_{j+2}\}$. The base method, ignoring the optimal selection of $\lambda^2$, reconstructs $u_{j+1/2}^{-}$ and $u_{j+1/2}^{+}$ using a three-point reconstruction stencil. Hence, in order to compute the numerical flux at the cell boundary $x_{j+1/2}$, we require a total of four cell averages: $\{\bar{u}_{j-1},\bar{u}_j,\bar{u}_{j+1}, \bar{u}_{j+2}\}$; however, we show, below, that this same set of points can be used to approximate the optimal value of $\lambda^2$, which promotes the approximation from third to fourth order accuracy. While it may appear that the method requires an additional point, the effective stencil is identical to the one generated by the third order WENO-JS scheme. In other words, no additional points are required beyond what is already needed for the analogous classical scheme.

If we represent the cell average values using Taylor expansion, then it follows that
we can obtain the linear combination of cell averages $\bar{u}_k$ which approximate $u'''(x_{j+\frac12})$ and $u'(x_{j+\frac12})$, i.e., the derivatives of the function at the cell boundary $x_{j+\frac12}$, as follows:
\begin{align*}
     u'(x_{j+\frac12})  &=    \frac{1}{\Delta x} \left( \frac{1}{12}\bar{u}_{j-1} -\frac{5}{4}\bar{u}_{j} +\frac{5}{4}\bar{u}_{j+1} -\frac{1}{12}\bar{u}_{j+2} \right)  + \mathcal{O}(\Delta x^3) ,\\
     u'''(x_{j+\frac12}) &=  \frac{1}{\Delta x^3} \left( -\bar{u}_{j-1} + 3\bar{u}_{j} -3\bar{u}_{j+1} +\bar{u}_{j+2} \right)  + \mathcal{O}(\Delta x).
\end{align*}
We remark that because the shape parameter is used to promote the accuracy of the reconstructions in smooth regions, these difference approximations for the derivative do not need to account for wind-direction. Using these approximations, can compute the shape parameter
\begin{align*}
     \lambda^2 
     &= - \frac{ -\bar{u}_{j-1} + 3\bar{u}_{j} -3\bar{u}_{j+1} +\bar{u}_{j+2}}{ {\Delta x^2} \left( \bar{u}_{j-1} -{15}\bar{u}_{j} +{15}\bar{u}_{j+1} -\bar{u}_{j+2} \right)}, \\
     &= -\frac {u'''(x_{j+\frac12})}{12u'(x_{j+\frac12})} + \mathcal{O}(\Delta x),
\end{align*}
which meets the convergence criteria, shown in equation \eqref{eq:lam4rbf3}, required to promote the reconstruction to fourth order accuracy. In order to prevent a division by zero, in regions where the function data is ``flat", we instead use $$ u'(x_{j+\frac12}) \approx \frac{1}{\Delta x} \left( \frac{1}{12}\bar{u}_{j-1} -\frac{5}{4}\bar{u}_{j} +\frac{5}{4}\bar{u}_{j+1} -\frac{1}{12}\bar{u}_{j+2} \right)  +  \varepsilon, \quad \varepsilon: = \varepsilon(\Delta x), $$
with $\varepsilon$ having the same sign as the first group of terms involving differences of the cell average data.

\section{Stencil Coefficients for the WENO-RBF Methods} \label{app:RBF_coeffs}

This section provides the stencil coefficients used by the WENO-RBF schemes presented in this work, which assume the basis function is a Gaussian (see section \ref{sec: Superconvergent schemes}). We treat the shape parameter(s), i.e., $\lambda^2$, as input in the reconstruction procedure. For simplicity, we present the expressions using a single shape parameter, but, in general, each (sub)stencil may be associated with its own shape parameter. An example calculation, for the shape parameter, in the fixed-stencil RBF3 scheme, is provided in Appendix \ref{app:lambda}.

\subsection{Three-Point Scheme}
\label{app_sub:three-point}
 
 Using the global stencil of cell-averages $S_{3} := \{ \bar{u}_{j-1}, \bar{u}_{j}, \bar{u}_{j+1} \}$ one obtains, following the procedure in section \ref{sec: Superconvergent schemes}, the fixed stencil reconstruction of the form $$ u_{j+1/2}^{S_3} = C_{-1} \bar{u}_{j-1} + C_{0} \bar{u}_{j} + C_{1}
 \bar{u}_{j+1}. $$ To improve the efficiency of the method and simplify the implementation of the method, we chose to Taylor expand the expressions coefficients (up to the order of the global stencil), which results in the following: 
 \begin{align*}
     C_{-1} &= -\frac{1}{6} - \frac{1}{3} \lambda^2 \Delta x^2 + \left( \lambda^2 \Delta x^2 \right)^{2} - \frac{5}{9} \left( \lambda^2 \Delta x^2 \right)^{3},\\
     C_{0}  &= \frac{5}{6} - \frac{1}{3} \lambda^2 \Delta x^2 + \frac{5}{6}\left( \lambda^2 \Delta x^2 \right)^{2} - \frac{5}{9} \left( \lambda^2 \Delta x^2 \right)^{3},\\
     C_{1}  &= \frac{1}{3} - \frac{2}{3} \lambda^2 \Delta x^2 + \frac{1}{6}\left( \lambda^2 \Delta x^2 \right)^{2} - \frac{5}{9} \left( \lambda^2 \Delta x^2 \right)^{3}. 
 \end{align*}
 
Similarly, using the two-point substencils $S_0 := \{ \bar{u}_{j-1}, \bar{u}_{j} \}$ and $S_1 := \{ \bar{u}_{j}, \bar{u}_{j+1} \}$, one obtains the reconstructions
\begin{align*}
    u_{j+1/2}^{(0)} &= c_{0}^{0} \bar{u}_{j-1} + c_{1}^{0} \bar{u}_{j}, \\
    u_{j+1/2}^{(1)} &= c_{0}^{1} \bar{u}_{j} + c_{1}^{1} \bar{u}_{j+1},
\end{align*}
with the corresponding (expanded) coefficients given by
\begin{align*}
    c_{0}^{0} &= -\frac{1}{2} - \frac{2}{3} \left( \lambda^2 \Delta x^2 \right)^{2} + \left( \lambda^2 \Delta x^2 \right)^{3},\\
    c_{1}^{0} &= \frac{3}{2} - 2\lambda^2 \Delta x^2 + \left( \lambda^2 \Delta x^2 \right)^{2}, \\
    c_{0}^{1} &= \frac{1}{2} + \frac{1}{2} \lambda^2 \Delta x^2 - \frac{1}{12}\left( \lambda^2 \Delta x^2 \right)^{2} - \frac{1}{4}\left( \lambda^2 \Delta x^2 \right)^{3},\\
    c_{1}^{1} &= \frac{1}{2} + \frac{1}{2} \lambda^2 \Delta x^2 - \frac{1}{12}\left( \lambda^2 \Delta x^2 \right)^{2} - \frac{1}{4}\left( \lambda^2 \Delta x^2 \right)^{3}. 
\end{align*}

\subsection{Four-Point Scheme}
\label{app_sub:four-point}

 Using the global stencil of cell-averages $S_{4} := \{ \bar{u}_{j-1}, \cdots, \bar{u}_{j+2} \}$ one obtains, following the procedure in section \ref{sec: Superconvergent schemes}, the fixed stencil reconstruction of the form $$ u_{j+1/2}^{S_4} = C_{-1} \bar{u}_{j-1} + C_{0} \bar{u}_{j} + C_{1}
 \bar{u}_{j+1} + C_{2} \bar{u}_{j+2}. $$ Performing Taylor expansions on the set coefficients, as in Appendix \ref{app_sub:three-point}, we obtain 
 \begin{align*}
     C_{-1} &= -\frac{1}{12} - \frac{1}{3} \lambda^2 \Delta x^2 - \frac{1}{3} \left( \lambda^2 \Delta x^2 \right)^{2} + \frac{4}{9} \left( \lambda^2 \Delta x^2 \right)^{4},\\
     C_{0}  &= \frac{7}{12} + \frac{1}{3} \lambda^2 \Delta x^2 - \frac{2}{3}\left( \lambda^2 \Delta x^2 \right)^{2} - \frac{1}{9} \left( \lambda^2 \Delta x^2 \right)^{4},\\
     C_{1}  &= \frac{7}{12} + \frac{1}{3} \lambda^2 \Delta x^2 - \frac{2}{3}\left( \lambda^2 \Delta x^2 \right)^{2} - \frac{1}{9} \left( \lambda^2 \Delta x^2 \right)^{4}, \\
     C_{2}  &= -\frac{1}{12} - \frac{1}{3} \lambda^2 \Delta x^2 - \frac{1}{3} \left( \lambda^2 \Delta x^2 \right)^{2} + \frac{4}{9} \left( \lambda^2 \Delta x^2 \right)^{4}
 \end{align*}
 
On each of the two-point substencils $S_0 := \{ \bar{u}_{j-1}, \bar{u}_{j} \}$, $S_1 := \{ \bar{u}_{j}, \bar{u}_{j+1} \}$, and $S_2 := \{ \bar{u}_{j+1}, \bar{u}_{j+2} \}$, one obtains the reconstructions
\begin{align*}
    u_{j+1/2}^{(0)} &= c_{0}^{0} \bar{u}_{j-1} + c_{1}^{0} \bar{u}_{j}, \\
    u_{j+1/2}^{(1)} &= c_{0}^{1} \bar{u}_{j} + c_{1}^{1} \bar{u}_{j+1}, \\
    u_{j+1/2}^{(2)} &= c_{0}^{2} \bar{u}_{j+1} + c_{1}^{2} \bar{u}_{j+2}
\end{align*}
with the corresponding (expanded) coefficients given by
\begin{align*}
    c_{0}^{0} &= -\frac{1}{2} - \frac{2}{3} \left( \lambda^2 \Delta x^2 \right)^{2} + \left( \lambda^2 \Delta x^2 \right)^{3} - \frac{103}{180} \left( \lambda^2 \Delta x^2 \right)^{4} ,\\
    c_{1}^{0} &= \frac{3}{2} - 2\lambda^2 \Delta x^2 + \left( \lambda^2 \Delta x^2 \right)^{2} - \frac{1}{5} \left( \lambda^2 \Delta x^2 \right)^{4}, \\
    c_{0}^{1} &= \frac{1}{2} + \frac{1}{2} \lambda^2 \Delta x^2 - \frac{1}{12}\left( \lambda^2 \Delta x^2 \right)^{2} - \frac{1}{4}\left( \lambda^2 \Delta x^2 \right)^{3} + \frac{7}{720} \left( \lambda^2 \Delta x^2 \right)^{4},\\
    c_{1}^{1} &= \frac{1}{2} + \frac{1}{2} \lambda^2 \Delta x^2 - \frac{1}{12}\left( \lambda^2 \Delta x^2 \right)^{2} - \frac{1}{4}\left( \lambda^2 \Delta x^2 \right)^{3} + \frac{7}{720} \left( \lambda^2 \Delta x^2 \right)^{4}, \\
    c_{0}^{2} &= \frac{3}{2} - 2\lambda^2 \Delta x^2 + 2\left( \lambda^2 \Delta x^2 \right)^{2} - \left( \lambda^2 \Delta x^2 \right)^{3} + \frac{13}{60} \left( \lambda^2 \Delta x^2 \right)^{4},\\
    c_{1}^{2} &= -\frac{1}{2} + \frac{1}{3}\left( \lambda^2 \Delta x^2 \right)^{2} - \frac{7}{45} \left( \lambda^2 \Delta x^2 \right)^{4}.
\end{align*}

\end{appendices}

\end{document}